\documentclass[12pt, reqno]{amsart}
\usepackage{amsmath, amsthm, amscd, amsfonts, amssymb, graphicx, color}
\usepackage[bookmarksnumbered, colorlinks, plainpages]{hyperref}
\hypersetup{colorlinks=true,linkcolor=red, anchorcolor=green, citecolor=cyan, urlcolor=red, filecolor=magenta, pdftoolbar=true}

\usepackage{amsmath,amssymb,colordvi,mathrsfs}
\usepackage{verbatim}
\usepackage{color}
\usepackage{tikz}
\newtheorem{theorem}{Theorem}
\newtheorem{corollary}[theorem]{Corollary}
\newtheorem{lemma}[theorem]{Lemma}
\newtheorem{example}[theorem]{Example}
\newtheorem{proposition}[theorem]{Proposition}
\newtheorem{definition}[theorem]{Definition}

\newtheorem{remark}[theorem]{\it Remark}

\newcommand{\CaixaPreta}{\vrule Depth0pt height5pt width5pt}
\newcommand{\bgproof}{\noindent {\bf Proof.} \hspace{2mm}}
\newcommand{\edproof}{\hfill \CaixaPreta \vspace{3mm}}

\def\NN{\mathbb N}
\def\CC{\mathbb C}

\def\R{\mathbb R}
\def\RR{\mathbb R}

\begin{document}

\title[Sharp extensions and algebraic properties]{Sharp extensions and algebraic properties for solution families of vector-valued differential equations}

\author[L. Abadias, C. Lizama, P. J. Miana]{Luciano Abadias$^1$, Carlos Lizama$^2$$^{*}$ and Pedro J. Miana$^1$}

\address{$^{1}$Departamento de Matem\'aticas, Instituto Universitario de Matem\'aticas y Aplicaciones, Universidad de Zaragoza, 50009 Zaragoza, Spain;
\newline
 L. Abadias and P. J. Miana  have been partially supported by Project MTM2013-42105-P, DGI-FEDER, of the MCYTS; Project E-64, D.G. Arag\'on, and UZCUD2014-CIE-09, Universidad de Zaragoza, Spain.}
\email{\textcolor[rgb]{0.00,0.00,0.84}{labadias@unizar.es;
pjmiana@unizar.es}}

\address{$^{2}$ Departamento de Matem\'atica y Ciencia de la Computaci\'on, Universidad de Santiago de Chile, Casilla 307-Correo 2, Santiago-Chile, Chile;
\newline
C. Lizama has been partially supported by DICYT, Universidad de Santiago de Chile; Project
CONICYT-PIA ACT1112 Stochastic Analysis Research Network; FONDECYT 1140258 and Ministerio de Educaci\'{o}n CEI Iberus (Spain).}
\email{\textcolor[rgb]{0.00,0.00,0.84}{carlos.lizama@usach.cl}}




\date{Received: xxxxxx; Revised: yyyyyy; Accepted: zzzzzz.
\newline \indent $^{*}$ Corresponding author}

\begin{abstract}
In this paper we show the unexpected property that extension from local to global without loss of regularity holds for the solutions of a wide class of vector-valued differential equations, in particular for the class of fractional abstract Cauchy problems in the subdiffusive case.  The main technique is the use of the algebraic structure of these solutions, which are defined by new versions of functional equations defining solution families of bounded operators. The convolution product and the double Laplace transform for functions of two variables are useful tools which we apply also to extend these solutions. Finally we illustrate our results with different concrete examples.
\end{abstract} \maketitle

\date{}

\maketitle

\section{Introduction}
\setcounter{theorem}{0}

Let $A$ be a closed linear operator with domain $D(A)$ defined in a complex Banach space $X$ and $0<\tau\leq \infty.$ Suppose that $A$ is the generator of a {\it local} $C_0$-semigroup $\{T(t)\}_{t \in [0,\tau)}$ or, equivalently,  the first order Cauchy problem
\begin{equation}\label{mainequation}
\left\{\begin{array}{ll}
u'(t)=Au(t)+x,&0\leq t < \tau \\
u(0)=0&
\end{array} \right.
\end{equation}
has a unique solution $u\in C^1([0,\tau),X)\cap C([0,\tau), D(A)),$ i.e. is {\it locally well-posed.} Then it is well known that $A$ is the generator of a {\it global} $C_0$-semigroup $\{T(t)\}_{t \geq 0}$, i.e. the problem is {\it globally well-posed,} see  \cite[Theorem 1.2]{Arendt-Elmen-Ke}, \cite[Section 3.1]{ABHN} and also \cite{Wa-Gao}. We observe that this dynamic behavior of the solution for the Cauchy problem \eqref{mainequation}, i.e. the extension property from local to global without loss of regularity, heavily depends on the {\it translation in time property} of the Cauchy's functional equation, namely $T(t+s)=T(t)T(s), \,\, t,s\geq 0.$
In contrast, this extension property is not more true for the class of local integrated semigroups  (\cite[Example 4.6]{Arendt-Elmen-Ke}). {Furthermore,} we have that if $A$ generates a local $(1*k)$-convoluted semigroup on $[0,\tau)$, then $A$ generates a local $(1*k^{*n})$-convoluted semigroup  on an interval $[0,n\tau)$, see \cite[Theorem 4.4]{Ke-Mi-La14} (\cite[Theorem 3.3]{Mi-Po14}). In other words, in these cases there is evolution with jumps of regularity and naturally the need of regularize the family of operators appears (in the sense of convolution) in order to have extension. { In all of these cases,} the property of translation in time of the associated functional equation is strongly connected with the problem of extension from a short to a long interval of definition for the corresponding family of operators.

Fractional diffusion equations are widely used to describe anomalous diffusion processes. From the point of view of operator theoretical-methods for Partial Differential Equations, subdiffusion phenomena is modeled naturally by means of fractional Cauchy problems in the form
\begin{equation}\label{eqmain2}
\left\{\begin{array}{ll}
D^{\alpha}_t u(t)=Au(t)+x,&0\leq t < \tau, \\
u(0)=0,&
\end{array} \right.
\end{equation}
where $0<\alpha < 1$ and the fractional derivative is taken in the Caputo sense, see \cite{Bajlekova, Po99}. In \cite{Li} the existence of solutions of fractional Cauchy problems is studied in detail, and  in the reference \cite{LiPe} for the superdiffusive case  $1<\alpha<2$. Suppose that $A$ generates a local one-parameter family of bounded operators that makes the equation \eqref{eqmain2} locally well posed. The natural question that arises is:
Can be \eqref{eqmain2} globally well posed?.


We point out that \eqref{eqmain2} is included in the more general  Volterra type equation
\begin{equation}\label{eq1.1bb}
  u(t)= k(t)x + A \int_0^t a(t-s)u(s)ds, \quad t\in (0,\tau),
  \end{equation}
for the special choice of kernel $a(t)=k(t)= \frac{t^{\alpha-1}}{\Gamma(\alpha)},$ for $\alpha, t>0$ and where $\Gamma$ is the Euler Gamma function. { This Volterra type equation in case $k(t)\equiv 1$ has been deeply treated in the monograph \cite{Pr93} by J. Pr\"uss. Further relevant studies have been done in the monographs by M. Kosti\'c \cite{Kostic1, Kostic2}. See also the references \cite{Ke-Li-Mi09, Ke-Mi-La14, Kostic, Li-Sun1, Liz-Po12, Liz, Mei-Peng1} and \cite{Mi-Po14}  for related work.} Therefore, we can set  our problem in a more general context: Classify the classes of pairs $(a,k)$ where extension of \eqref{eq1.1bb} from local to global  and without loss of regularity holds.

We note that under certain conditions on the scalar-valued kernels $a$ and $k$,  {\it  well posedness} of the Volterra equation \eqref{eq1.1bb} is equivalent to the existence of a one-parameter and strongly continuous family of bounded operators $\{S(t)\}_{t > 0}$ that satisfies a functional equation in the form
\begin{equation}
\begin{array}{lcl}\label{eq1.1}
S(s)\displaystyle \int_0^t a(t-\tau)S(\tau)d\tau &-& \displaystyle  S(t) \int_0^s a(s-\tau)S(\tau)d\tau  \\
=k(s)\displaystyle \int_0^t a(t-\tau)S(\tau)d\tau &-& k(t)\displaystyle \int_0^s a(s-\tau)S(\tau)d\tau,
\end{array}
\end{equation}
for $t,s > 0,$ see  \cite[Theorem 3.1]{Liz-Po12}.

\bigskip
Explicitly, in this paper we will study the following  questions:\\

{\bf (Q1)} (Evolution  with or without jumps of regularity):  If $A$ is the generator  of a local regularized resolvent family on the interval $(0, \tau)$, is $A$ also the generator of an local extended regularized resolvent family on the interval $(0, (n+1)\tau)$ for $n\in \NN$; in particular for which class of pairs of kernels $(a,k)$ we have that  if $A$ is the generator of a local $(a,k)$-regularized resolvent family, then $A$ is  the generator of a global $(a,k)$-regularized resolvent family? (Sections \ref{local} and \ref{jump}).\\

{\bf (Q2)} (Time translation): Determine the class of pairs $(a,k)$ for which is possible to find an equivalent  and explicit expression for equation (\ref{eq1.1}) in terms of the sum $t+s$ instead of $t$ and $s$. (Section \ref{algebraic}) \\

We will first answer globally the problem of evolution  {\bf (Q1)}, that is, the possibility to extend the family of operators $S(t)$ from the interval $(0,\tau)$ to the whole semiaxis $(0,\infty).$ More precisely, we prove: If $A$ is the generator of a local $(a,k)$-regularized resolvent family on $(0,\tau),$ then $A$ is the generator of a local $(a, (a*k)^{*n} *k)$-regularized resolvent family in $(0,(n+1)\tau),$ see Theorem \ref{local1}. We remark that this problem, which has been studied in a series of papers in the last years, is settled here in a simple way, making transparent the process of regularization needed in each step of the extension. In particular, our result intersects the papers \cite{Arendt-Elmen-Ke}, \cite{Ke-Mi-La14} and \cite{Mi-Po14} where the problem of extension for local integrated semigroups, local convoluted semigroups and  local convoluted cosine functions is studied, respectively. { The question about extension for local convoluted semigroups and local convoluted cosine funtions, which is resolved in \cite{Ke-Mi-La14} and \cite{Mi-Po14} respectively, had been cited previously in the paragraph directly preceding \cite[Theorem 1.2.7]{Kostic1}. Results related to the extension of local $C$-regularized semigroups and local $C$-regularized cosine functions appeared by the the first time in \cite{Wa-Gao}.}

However, note that if we restrict this result to the  $C_0$-semigroup family, we do not obtain  evolution conserving regularity. Under some conditions on $a$ and $k$, we  improve this result and obtain extension without jump { on} regularity in  Theorem \ref{th-local1}. We use this Theorem to prove one of the main results in  this paper concerning the fractional equation \eqref{eqmain2}: If $A$ is the generator of a local $(g_{\alpha},g_{\alpha})$-regularized resolvent family on $[0,\tau),$  with $0<\alpha<1,$ then $A$ is the generator of a global $(g_{\alpha}, g_{\alpha})$-regularized resolvent family in $[0,\infty),$ see Corollary \ref{co-local1}.
These results recover and widely extend  the property of evolution without jumps of regularity from the case of the solutions of first order Cauchy problems to the case of fractional subdiffusive models.

In this paper, we are able to completely solve {\bf (Q2)} establishing an equivalent functional equation to (\ref{eq1.1}), which defines global $(a,k)$-regularized resolvent families,  in the following form:
\begin{equation}\label{eqf}
\begin{array}{r}
 \displaystyle \int_t^{t+s} \int_0^r k(s+t-r)a(r-\tau)S(\tau)d\tau dr - \displaystyle \int_0^s \int_0^r k(s+t-r)a(r-\tau)S(\tau)d\tau dr \\ \\ = \displaystyle \int_0^t \int_0^s a(t+s-r_1-r_2)S(r_1)S(r_2)dr_1 dr_2,
\end{array}
\end{equation}
for $t,s\ge 0$ (Theorem \ref{th4.1}). The above formula widely solves  the  problem of time translation, not only extending all the several results existing in the literature, but also proposing and finding a better expression for older and new  cases. For example, we will see that for $a(t)={\frac{t^{\alpha-1}}{\Gamma(\alpha)}}$ and $k(t)= 1,$ for $t>0$, the border cases $\alpha=1$ and $\alpha=2$ are naturally included in our functional formula, unifying the cases $0<\alpha<1$ and $1<\alpha<2$ mentioned before. We point out that  functional equation (\ref{eqf})  inspired the way to define extensions in local $(a,k)$-regularized resolvent family commented above.

In last years, special interest has appeared in the study of algebraic function equations {\bf (Q2)}  for $(a,k)$-regularized resolvent family only with $a(t)=\frac{t^{\alpha-1}}{\Gamma(\alpha)}$  for $t>0$ (due mainly for its connection with fractional differential equations). First results in this line appeared in \cite{Peng-Li-12} where an equivalent functional equation to (\ref{eq1.1})  in case $a(t)=\frac{t^{\alpha-1}}{\Gamma(\alpha)}$ for $0<\alpha<1,$ and $k(t)=1,$  $(t>0)$ is given in \cite[ Formula (2.1)]{Peng-Li-12}.
After that,  a similar result for the case $a(t)=k(t)= \frac{t^{\alpha-1}}{\Gamma(\alpha)},$
 is shown in \cite[Proposition 2.2]{Mei-Pe-Zh13}.
  In the recent paper \cite{Li-Sun1}, a further extension of  known result   in case $a(t)=\frac{t^{\alpha-1}}{\Gamma(\alpha)}$ and $k(t)= \frac{t^{\beta}}{\Gamma(\beta+1)}$ was proved (\cite[Theorem 5]{Li-Sun1}).
   We note that some restriction should be imposed for $\alpha\geq 1,$ see Example \ref{ex4.9}. A further generalization, this time in case $a(t)=\frac{t^{\alpha-1}}{\Gamma(\alpha)}$ for $0<\alpha<1$ and $k(t)= \int_0^t K(s)ds$ was successfully obtained  in \cite[Theorem 8]{Mei-Peng1}.
 Finally we point out that in a very recent work \cite{Mei-Peng-Jia}, the authors have discovered a functional equation which cover the case $a(t)= \frac{t^{\alpha-1}}{\Gamma(\alpha)}$ and $k(t)= 1 $ ($t>0$) in the super diffusive case $1<\alpha<2$ (\cite[Definition 3.1]{Mei-Peng-Jia}).


To handle both questions {\bf (Q1, Q2)}, we introduce an original technique in this context: we consider scalar and vector-valued functions of two variables. In second section we work with convolution product $*_2$ (see formula (\ref{de2.2})), and we prove some needed technical results that play a key role in the paper. In third section, we see results about simple and double Laplace transform, properties that this transforms verifies and in which appear the above convolutions products (see, for example, Theorem \ref{le3.4}).

Double Laplace transform is an efficient and known tool to solve scalar differential equations in two variables, see for example \cite[Chapitre IV.15.3]{DP},  \cite{Jaeger}, \cite[pp 226--228]{Sne}. Other interesting applications of double Laplace transform is to supply integral formulae (\cite[Chapitre IV.15.1]{DP}) and bilinear expansions (\cite[Chapitre IV.15.2]{DP}). In \cite{Str} the structure of closed ideals of convolution algebra $L^1(\RR^n)$ is studied and the Laplace transform for functions of several variables is also considered.

Finally, in the last section we illustrate our main results with some particular examples: we considerer some  $(a, k)$-regularized resolvent families, some of which are related with the solution of fractional Cauchy problems to conclude that the extension process is possible with or without jump of regularity. We also give new functional equations obtained as a consequence of the formula (\ref{eqf}) for known one-parametric family of operators: $C_0$-semigroups, cosine families, convoluted semigroups and resolvent families.

\section{Convolution products in one and two variables}
\setcounter{theorem}{0}
\setcounter{equation}{0}

In this section we state some technical results for convolution products (in one and two variables) that we use to prove the relevant results in the following sections. The convolution product in several variables have been considered in some relevant fields in Mathematical Analysis (see for example \cite{Str}). However, the convolution product in two variables is a new tool to apply to $(a,k)$-regularized resolvent families.

We denote $\R_+=[0,+\infty)$; $\R=(-\infty, +\infty)$, $\R_+^2=\R_+\times\R_+$ and $\R^2=\R\times\R.$  We consider the space of locally integrable functions  in one and two variables, $L^1_{loc}(\R_+),$ and $L^1_{loc}(\R_+^2)$. The space $\mathcal{C}^n(\mathbb{R}_+)$ is formed with continuous functions $f:\R_+\to \CC$  such $f^{(j)}$ is continuous for $0\le j\le n$ for $n\ge 0$. { Some of the above spaces will also be considered in $(0,\tau)$ or $[0,\tau)$ instead $\R_+,$ with $\tau>0.$}

Let $f,g:\R_+\to\CC$, we write $f_t(s):=f(s+t)\chi_{[-s,+\infty)}(t)$ for $t\in\R;$ $f^+:\R_+^2\to \CC$ the function given by $f^+(t,s):=f(t+s);$ $f^-:\R^2\to \CC$ the function given by $f^-(t,s):=f(|t-s|);$ $f\otimes g:\R_+^2\to \CC$ by $f\otimes g(t,s):= f(t)g(s)$ for $(t,s)\in \R_+^2$  and \begin{equation}\label{convo}f*g(t)=\int_0^{t}f(t-s)g(s)\,ds, \qquad t>0,\end{equation} the usual convolution product, for the functions $f,g$ where the product is convergent. We write $f^{*2}$ instead $f*f$ and then $f^{*n}=f*(f^{*(n-1)})$ for $n\geq 2$ is the $n$-fold convolution power of $f.$

For $F,G:\R^2_+\to\CC$ be given, we define the  convolution product in two variables by \begin{equation} \label{de2.2}
F*_2G(t,s):=\int_0^t\int_0^s F(t-u,s-v)G(u,v)\,dv\,du, \qquad t,s>0.
\end{equation}
whenever is well defined. This product is commutative and associative, see \cite[Formula (13.9)]{DP} and \cite[Formula (3-18-19)]{Sne}.

  We define functions $g_{\alpha}(t):=\frac{t^{\alpha-1}}{\Gamma(\alpha)}$,  $e_\lambda(t):=e^{-\lambda t}$ and $e_{\lambda, \mu}(t,s):=e^{-\lambda t-\mu s}=e_{\lambda}\otimes e_{\mu}(t,s)$ for  $\alpha\in \R\backslash\{0,-1,-2,-3\dots\}$, $\lambda, \mu\in \CC $ and $t,s>0$.  Note that $(e_\lambda)^+=e_{\lambda, \lambda}$ for $\lambda\in \CC.$ It is direct to check the following well-known identities:
  \begin{eqnarray*}
  g_\alpha\ast g_\beta&=&g_{\alpha+\beta}, \qquad \alpha, \beta>0;\cr
  e_\lambda \ast e_\mu=&=&{1\over \lambda-\mu}\left(e_\mu-e_\lambda\right), \qquad \lambda\not =\mu;\cr
  e_{\lambda,\lambda'} \ast_2 e_{\mu,\mu'}=&=&{1\over (\lambda-\mu)(\lambda'-\mu')}\left(e_{\lambda,\lambda'}-e_{\lambda, \mu'}-e_{\mu,\lambda'}+e_{\mu, \mu'}\right), \quad \lambda\not =\mu,\ \lambda'\not =\mu'.
  \end{eqnarray*}

 The way that $\ast$ and $\ast_2$ interact with operators $\otimes$, $(\cdot)_t$ and $(\cdot)^+$ is shown in the next theorem.

\begin{theorem}\label{le2.2} Let $f,g,h,j\in L^1_{loc}(\R_+).$ Then
\begin{itemize}
\item[(i)] $(f\otimes g)\ast_2(h\otimes j)= (f\ast h)\otimes(g\ast j)$.
\item[(ii)]
$(g^+*_2(f\otimes h))(t,s)=h*(f*g)_t(s)-f_t*(h*g)(s),$  for $t,s\geq 0.
$
\item[(iii)] \small{$$(f^+\ast_2 g^+)(t,s)=\left\{\begin{array}{ll}
(f_t\ast {\mathcal M}(g))(s)+s(f_s\ast g_s)(t-s)+ { ({\mathcal M}(f)\ast g_t)(s), } &0\leq s \le t, \\
(f_s\ast {\mathcal M}(g))(t)+t(f_t\ast g_t)(s-t)+ {({\mathcal M}(f)\ast g_s)(t), }&0\leq t \le s,
\end{array} \right.
$$}
where $ {\mathcal M}(g)(s):=sg(s)$ for $s\in \R_+$.
\end{itemize}
\end{theorem}
\bgproof The proof of  part (i) is straightforward. To show (ii), note that if $g\in L^1_{loc}(\R_+)$ then $g^+\in L^1_{loc}(\R_+^2).$ We  change variables to obtain the following equalities: \small{\begin{displaymath}\begin{array}{l}
\displaystyle\int_0^t\int_0^s g(t+s-r_1-r_2)f(r_1)h(r_2)\,dr_1\,dr_2=\int_0^t\int_0^s g(v+z)f(t-v)h(s-z)\,dv\,dz \\ \\
\displaystyle=\int_0^s h(s-z)\int_z^{t+z}f(t+z-u)g(u)\,du\,dz \\ \\
\displaystyle=\int_0^s h(s-z)\int_0^{t+z}f(t+z-u)g(u)\,du\,dz -\int_0^s h(s-z)\int_0^{z}f(t+z-u)g(u)\,du\,dz .
\end{array}\end{displaymath}}
Now, we apply Fubini theorem and  change of variable $s-z=r-u$ to get \small{\begin{displaymath}\begin{array}{l}
\displaystyle\int_0^s h(s-z)\int_0^{z}f(t+z-u)g(u)\,du\,dz=\int_0^s g(u)\int_u^{s}f(t+z-u)h(s-z)\,dz\,du \\ \\
\displaystyle=\int_0^s g(u)\int_u^{s}f(t+s-r)h(r-u)\,dr\,du =\int_0^s f(t+s-r)\int_0^{r}h(r-u)g(u)\,du\,dr .
\end{array}\end{displaymath}}
We express the above integrals in terms of convolution products to conclude the claim. The proof of (iii) is similar to the proof of part (ii).
\edproof

\begin{corollary}\label{cor2.3} Let $f,g,h\in L^1_{loc}(\R_+).$ Then
\begin{itemize}
\item[(i)] $(g^+*_2(f\otimes h))(t,s)=f*(h*g)_s(t)-h_s*(f*g)(t),$ for $t,s\geq 0.$
\item[(ii)]  $(g^+ *_2(f\otimes f))(t,s)=f*(f*g)_t(s)-f_t*(f*g)(s),$ for $t,s\geq 0.$
\end{itemize}
\end{corollary}

\bgproof (i) We apply the identity
$(g^+*_2(f\otimes h))(t,s)=(g^+*_2(h\otimes f))(s,t)
$ for $(t,s)\in \R_+^2$ and  the Theorem \ref{le2.2} (ii).
\edproof

Next lemma extends \cite[Lemma 2.1]{Ke-Mi-La14} and will be applied several times in this paper.

\begin{lemma}\label{le5.1} Take $0\leq \tau \leq t$ and $f,g,h\in L^1_{loc}(\R_+).$ Then
\begin{eqnarray*}
 \int_0^{t-\tau}h(t-s)(g*f)(s)\,ds &+& \int_0^{\tau}f(t-s)(g*h)(s)\,ds \\ &=& (f*g*h)(t)- g^+\ast_2(f\otimes h)(t-\tau,\tau).
 \end{eqnarray*}
\end{lemma}
\bgproof We use Fubini's theorem and change of variables to obtain
$$
\begin{array}{l}
 \displaystyle (f*g*h)(t)-\int_0^{\tau}f(t-r)\int_0^rg(r-s)h(s)\,ds\,dr \\ \\
\displaystyle =\int_0^{t}f(r)\int_0^{t-r}g(t-r-s)h(s)\,ds\,dr  - \int_{t-\tau}^{t}f(r)\int_0^{t-r}g(t-r-s)h(s)\,ds\,dr \\ \\
\displaystyle =\int_0^{t-\tau}f(r)\int_0^{t-r}g(t-r-s)h(s)\,ds\,dr \\ \\
\displaystyle =\int_0^{\tau}h(s)\int_{0}^{t-\tau} g(t-s-r)f(r)\,dr\,ds + \int_{\tau}^{t}h(s)\int_{0}^{t-s} g(t-s-r)f(r)\,dr\,ds \\ \\
\displaystyle =\int_0^{\tau}\int_0^{t-\tau}g(t-s-r)h(s)f(r)\,dr\,ds  + \int_0^{t-\tau}h(t-s)(g*f)(s)\,ds,
\end{array}
$$
for $t\in\R_+.${ This} proves the claim.
\edproof

\begin{remark}{\rm
Let $X$ be a Banach space and $$L^1_{loc}(\R_+,X):=\{f:\R_+\to X\,:\, f \text{ is Bochner integrable on }[0,\tau] \text{ for all } \tau>0\}.$$ We also consider $L^1_{loc}(\R_+^2,X)$ for functions defined in two variables.  The definitions of $\ast$ and $\ast_2$, (see (\ref{convo}) and (\ref{de2.2})), Theorem \ref{le2.2},  Corollary \ref{cor2.3} and Lemma \ref{le5.1} hold in the case that one function is vector valued  into $X.$ The proof of these analogous results { involves the ideas already employed in the scalar case.}}
\end{remark}

\section{Laplace transform in one and two variables}
\setcounter{theorem}{0}
\setcounter{equation}{0}

In this section we study properties concerned with the Laplace transform for functions in the above spaces. We say that $f\in L^1_{loc}(\R_+,X)$ is a  Laplace transformable function if there exists $\omega_f\in\R$ such that  the usual Laplace transform $$\hat{f}(\lambda):=\int_0^{\infty}e^{-\lambda t}f(t)\,dt=\displaystyle\lim_{\tau\to\infty}\int_0^{\tau}e^{-\lambda t}f(t)\,dt, \qquad \Re \lambda>\omega_f,$$
is well-defined, see for example \cite[Section 1.4]{ABHN}. Let $f: \R^+\to X$ be absolutely continuous and differentiable a.e. { Note that in the scalar case $X=\mathbb{C}$ any absolutely continuous function defined for $t\geq 0$ $(t>0)$ is differentiable a.e. $t\geq 0$ $(t>0)$ because the space $\mathbb{C}$ has the Radon Nykodim property.} If $\Re \lambda>0$ and  $\hat{f'}(\lambda)$ exists then $\hat{f}(\lambda)$ exists and
\begin{equation}\label{derivada}
\hat{f'}(\lambda)=\lambda \hat {f}(\lambda)-f(0),
 \end{equation}
see \cite[Corollary 1.6.6]{ABHN}.

Similarly, we say that $F\in L^1_{loc}(\R_+^2,X)$ is a double Laplace transformable function (or 2-Laplace transformable) if there exist $\omega_{1,F},  \omega_{2,F}\in\R$ such that
 $$\mathcal{L}_2(F)(\lambda,\mu):=\int_0^{\infty}\int_0^{\infty}e^{-\lambda t}e^{-\mu s}F(t,s)\,ds\,dt:=\displaystyle\lim_{\tau\to\infty}\int_0^{\tau}\int_0^{\tau}e^{-\lambda t}e^{-\mu s}F(t,s)\,ds\,dt$$  converges for $ \Re \lambda>\omega_{1, F}$ and $\Re \mu>\omega_{2, F}$, see \cite[Chapitre IV]{DP} and \cite[Section 3.18]{Sne} in the scalar case; the  Laplace transform $\mathcal{L}_2$ is commonly named the double Laplace transform.

For further use we establish the following Theorem where we include some known identities of Laplace transform and double Laplace transform.

\begin{theorem}\label{convolution2} Let $f\in L^1_{loc}(\R_+,X)$ and $g\in L^1_{loc}(\R_+)$ be Laplace transformable functions. Then  the following  identities { hold.}
\begin{itemize}

\item[(i)] $ \mathcal{L}_2(f^+)(\lambda,\mu)=\frac{1}{\mu-\lambda}(\hat{f}(\lambda)-\hat{f}(\mu))$ for $\Re\,\lambda,\Re\,\mu>\omega_f$ with $\lambda \neq \mu.$
\item[(ii)] $\mathcal{L}_2(f^-)(\lambda,\mu)=\frac{1}{\lambda+\mu}(\hat{f}(\lambda)+\hat{f}(\mu))$ for $\Re\lambda, \Re\mu>\omega_f$ with $ \Re(\lambda+\mu)>0.$

\item[(iii)] $\mathcal{L}_2(f\otimes g)(\lambda,\mu)= \hat{f}(\lambda)\hat{g}(\mu)$ for $\Re \lambda>\omega_f$ and $\Re \mu>\omega_g.$
\end{itemize}

Let
$F\in L^1_{loc}(\R_+^2,X)$ and $ G\in L^1_{loc}(\R^2_+)$ be double Laplace transformable functions. Then the following identity holds.
\begin{itemize}

\item[(iv)] $\mathcal{L}_2(F*_2G)(\lambda,\mu)=\mathcal{L}_2(G)(\lambda,\mu)\mathcal{L}_2(F)(\lambda,\mu),$ for $\Re\,\lambda>\max(\omega_{1,F}, \omega_{1, G})$ and  $\Re\,\mu>\max(\omega_{2,F}, \omega_{2, G})$.

\end{itemize}

\end{theorem}

\bgproof
The proof of (i) appears in \cite[pp 221-222]{Sne}; the proof of (ii) in \cite[pp 223-224]{Sne} and  the equality (iii) appears in \cite[(3-18-4)]{Sne}. Finally the equality (iv) is straightforward and it is commented in \cite[(3-18-20)]{Sne}.
\edproof


In what follows, given an  absolutely continuous and differentiable a.e. function { $c:(0,\infty)\to X$} we denote by $c'$ its derivative  and $c(0^+):=\displaystyle\lim_{t\to 0^+}c(t),$ whenever both limits exist.

\begin{theorem}\label{cor3.1} Let $c\in L^1_{loc}(\R_+,X)$ be an  absolutely continuous on  {$(0,\infty)$}, differentiable a.e and Laplace transformable function.
\begin{enumerate}
\item[(i)] If $(c')^+:\R_+^2\to X$ is 2-Laplace transformable, then
$$ \mathcal{L}_2((c')^+)(\lambda,\mu)=\frac{1}{\mu-\lambda}(\lambda\hat{c}(\lambda)-\mu\hat{c}(\mu)), \qquad \Re\,\lambda,\Re\,\mu>\omega_c, \qquad \lambda \neq \mu.$$

\item[(ii)] If $(c')^-:\R_+^2\to X$ is 2-Laplace transformable and $c(0^+)$ exists then $$\displaystyle \mathcal{L}_2((c')^-)(\lambda,\mu)=\frac{1}{\mu+\lambda}(\lambda\hat{c}(\lambda)+\mu\hat{c}(\mu)) - \frac{2c(0^+)}{\lambda + \mu}, \quad\Re\,\lambda, \Re\,\mu>\omega_c,\ \Re\,(\lambda+\mu)>0.$$
\end{enumerate}
\end{theorem}
\bgproof
(i) We integrate by parts to  obtain
\begin{displaymath}\begin{array}{l}
\displaystyle\int_0^{\infty}e^{-\lambda t}\int_0^{\infty}e^{-\mu s}c'(t+s)\,ds\,dt=\int_0^{\infty}e^{-(\lambda-\mu) t}\int_t^{\infty}e^{-\mu v}c'(v)\,dv\,dt\\
\displaystyle=\int_0^{\infty}e^{-(\lambda-\mu) t}\biggl( -c(t)e^{-\mu t}+\mu\int_t^{\infty}e^{-\mu v}c(v)\,dv \biggr)\,dt.\end{array}\end{displaymath}
We change the inner variable  to get the following equality,
$$
\int_0^{\infty}e^{-\lambda t}\int_0^{\infty}e^{-\mu s}c'(t+s)\,ds\,dt =
-\hat{c}(\lambda)+\mu\mathcal{L}_2(c^+)(\lambda,\mu)=\frac{\lambda\hat{c}(\lambda)-\mu\hat{c}(\mu)}{\mu-\lambda},
$$
for $\Re\,\lambda,\Re\,\mu>\omega_c,$ and $\lambda \neq \mu.$ (ii) For $\Re\,\lambda, \Re\,\mu>\omega_c,$ and  $\Re\,(\lambda+\mu)>0$, note that
\begin{displaymath}\begin{array}{l}
\displaystyle \mathcal{L}_2((c')^-)(\lambda,\mu) = \int_0^{\infty}e^{-\lambda t}\int_{-t}^{\infty}e^{-\mu (v+t)}(c')^{-}(t,v+t)\,dv\,dt  \\ \displaystyle= \int_0^{\infty}e^{-(\lambda+\mu) t}\int_0^{\infty}e^{-\mu v}c'(v)\,dv\,dt + \int_0^{\infty}e^{-(\lambda+\mu) t}\int_0^{t}e^{\mu v}c'(v)\,dv\,dt,
\end{array}\end{displaymath}
where we have changed the inner  variable. We integrate by parts to get that
$$
\int_0^{\infty}e^{-(\lambda+\mu) t}\int_0^{\infty}e^{-\mu v}c'(v)\,dv\,dt = \frac{1}{\lambda +\mu} \int_0^{\infty}e^{-\mu v}c'(v)\,dv= \frac{1}{\lambda +\mu}[-c(0^+)+ \mu \hat c(\mu)].
$$
On the other hand, we use  Fubini's theorem to obtain,
\begin{displaymath}\begin{array}{l}
\displaystyle\int_0^{\infty}e^{-(\lambda+\mu) t}\int_0^{t}e^{\mu v}c'(v)\,dv\,dt = \int_0^{\infty}e^{\mu v}c'(v)\int_v^{\infty}e^{-(\lambda+\mu) t}\,dt\,dv \\ \displaystyle= \frac{1}{\lambda +\mu} \int_0^{\infty}e^{-\lambda v}c'(v)\,dv = \frac{1}{\lambda +\mu}[-c(0^+)+ \lambda\hat c(\lambda)],
\end{array}\end{displaymath}
and we conclude the proof of the theorem.
\edproof

\begin{remark}{\rm In the case that the function $c'$ is a Laplace transformable function, the proof of Theorem \ref{cor3.1} is a straightforward consequence of Theorem \ref{convolution2} (i) and (ii) and the equality (\ref{derivada}). The interesting example $c= g_{1-\alpha}$ with $0<\alpha<1$ does not satisfy this condition, {however it  is absolutely continuous on $(0,\infty)$ and Laplace transformable and we can apply  Theorem \ref{cor3.1} (i) directly.}
}
\end{remark}

The following inversion theorem allows to express operators $(\quad)^+$ and $(\quad)^-$ in terms of   double convolution  products. { These equalities} play important roles in extension formulae which are obtained  in  next sections.

\begin{theorem}\label{le3.4} Let $a\in L^1_{loc}(\R_+)$ be a Laplace transformable function and suppose there exists $c\in L^1_{loc}(\R_+)$ { absolutely continuous on $(0,\infty)$} and Laplace transformable such that
\begin{equation}\label{Hypothesis1}
(a*c)(t)=1, \qquad t>0.
\end{equation}
\begin{enumerate}
\item[(i)] If $(c')^+$ is 2-Laplace transformable, then \small{$\displaystyle a^+=-((c')^+*_2 (a\otimes a)).$}
\item[(ii)] If $(c')^-$ is 2-Laplace transformable and $c(0^+)=0,$ then \small{$\displaystyle a^-=( (c')^-*_2 (a\otimes a)).$}
\end{enumerate}
\end{theorem}
\bgproof
(i) We use  Theorem \ref{convolution2} (i) and Theorem \ref{cor3.1} (i) to prove that:\begin{eqnarray*}
\mathcal{L}_2(a^+)(z,w)&=&\frac{\hat{a}(w)-\hat{a}(z)}{z-w}=\biggl(\frac{z\hat{c}(z)-w\hat{c}(w)}{z-w}\biggr)\frac{1}{zw\hat{c}(w)\hat{c}(z)}\\
&=&\biggl(\frac{z\hat{c}(z)-w\hat{c}(w)}{z-w}\biggr)\hat{a}(w)\hat{a}(z) =\mathcal{L}_2((c')^+*_2 (a\otimes a))(z,w).
\end{eqnarray*}
Due to the uniqueness of the Laplace transform (see for example \cite[p. 346]{DP}), we conclude the equality. The proof of part (ii) is similar and involves Theorem \ref{convolution2} (ii) and Theorem \ref{cor3.1} (ii).
\edproof

\begin{example}\label{co3.5}In the case that $c'$ is a Laplace transformable function, we apply Corollary \ref{cor2.3} (ii) to obtain an alternative proof of Theorem \ref{le3.4} (i). However the interesting example
$a=g_{\alpha}$ for $0<\alpha<1$  and  $c=g_{1-\alpha}$ does not satisfy this condition and the direct proof given in Theorem \ref{le3.4} (i) is needed. Note that $c'= g_{-\alpha}$ and the equality
 $\displaystyle g_{\alpha}^+=-((g_{-\alpha})^+ *_2 (g_{\alpha}\otimes g_{\alpha})), $ that is equivalent, after an algebraic manipulation, to the formula
 $$
 \frac{\alpha \sin{\alpha \pi}}{\pi} \int_0^t \int_0^s \frac{u^{\alpha-1}v^{\alpha-1}}{(t+u+s-v)^{\alpha +1}} ds\,dv = (t+s)^{\alpha -1}, \quad t,s >0.
 $$
\end{example}

\bigskip

Analogously, let $ S,\; T : \R_+\to {\mathcal B}(X) $ be strongly continuous operator families such that $S(\cdot)x,T(\cdot)x\in L^1_{loc}(\R_+,X),$ for any $x\in X.$ The operators $S, T$ are said Laplace-transformable functions if there { exists} $\omega\in \R$ such that  the Laplace transform of
$S$ {(respectively $T$)} $$ \hat S (\lambda)x = \displaystyle
\int_0^\infty e^{- \lambda t} S(t)x\, dt, \qquad\Re \lambda>\omega,$$ converges for $x\in X$, see for example \cite[Definition 3.1.4]{ABHN}. For $ h \in \mathbb{R}$ we shall
denote $ S_h $ the translation operator of $S$ given by  $ S_h (u) := S(u+ h) \chi_{[-h,
+\infty)}(u)$ for $u \in \mathbb{R}$ and the convolution product between $T$ and $S$, $T*S$, given by
$$
(T*S)(t)x := \int_0^t T(t-s)S(s)x\, ds, \qquad t>0, \quad x\in X.
$$
 If $g\in L^1_{loc}(\R_+)$ is a Laplace transformable scalar-valued function, then we define $g*S$ by $$(g*S)(t)x := \displaystyle \int_0^t g(t-s)S(s)x\, ds, \, t>0, \quad x\in X,$$ and $$
(g^+ *_2 (S\otimes S))(t,s)x := \int_0^t \int_0^s g(t+s-r_1-r_2)S(r_1)S(r_2)x\,dr_1\, dr_2, \quad t, s\geq 0,\ x\in X,
$$
 where $(S\otimes S)(t,s):=S(t)S(s)$ is the composition operator for $t, s\geq 0$. We recall  the following
identities given in  \cite[Lemma 4.1]{Ke-Li-Mi09}: for $
\lambda
> \mu
> \omega,$

\begin{equation}\label{eq2.1}
\hat S(\lambda) \hat T(\mu)x = \int_0^\infty e^{-\lambda t}
\int_0^\infty e^{-\mu s} S(t) T(s)x\, ds\, dt, \qquad x\in X,
\end{equation}
and
\begin{equation}\label{eq2.2}
 \frac{1}{ \mu - \lambda } ( \hat S(\lambda) - \hat S(\mu))x =
\int_0^{\infty} e^{- \lambda t} \int_0^\infty e^{-\mu s} S(t+s)x\, ds\,dt, \qquad x\in X.
\end{equation}
If $g: \R^+ \to \mathbb{C}$ is Laplace transformable, we also have
\begin{equation}\label{eq2.3}
\frac{1}{\mu - \lambda} \hat g(\mu)[ \hat S(\lambda) - \hat
S(\mu)]x = \int_0^\infty e^{- \lambda t} \int_0^\infty e^{-\mu s}
(g*S_t)(s)x\, ds\, dt, \qquad x\in X, \end{equation}
and
\begin{equation}\label{eq2.4}
\frac{1}{\mu - \lambda} \hat T(\mu)[ \hat g(\lambda) - \hat
g(\mu)]x = \int_0^\infty e^{- \lambda t} \int_0^\infty e^{-\mu s}
(T*g_t)(s)x\, ds\, dt, \qquad x\in X. \end{equation}
Defining $S(t)=S(-t)$ for $t<0,$ we have
\begin{equation}\label{eq2.6}
 \frac{1}{ \mu +\lambda } ( \hat S(\lambda) + \hat S(\mu))x =
\int_0^{\infty} e^{- \lambda t} \int_0^\infty e^{-\mu s} S(t-s)x\, ds\,dt, \qquad \lambda+\mu>0, \qquad x\in X.
\end{equation} In fact  equations \eqref{eq2.1}, \eqref{eq2.2}, \eqref{eq2.3} and \eqref{eq2.4} are valid for $\Re\, \lambda,\,\Re\,\mu>\omega$ with $\lambda\neq\mu,$ and \eqref{eq2.6} for $\Re\, \lambda,\Re\,\mu>\omega$ with $\Re(\lambda +\mu)>0.$

The following theorem shows how some double Laplace transforms of the double convolution product is also related  with the single Laplace transform.

\begin{proposition}\label{th3.1} Let $g\in L^1_{loc}(\R_+)$ and $S:\R_+\to \mathcal{B}(X)$ be a locally integrable and strongly continuous function, both Laplace transformable functions. Then the following identities hold \begin{enumerate}
\item[(i)]$
\displaystyle \mathcal{L}_2(g^+ *_2 (S\otimes S))(\lambda,\mu) = \frac{1}{\mu - \lambda }[\hat g(\lambda) -\hat g(\mu)]\hat S(\lambda) \hat S(\mu)$ for $\Re\,\lambda,\,\Re\,\mu>\omega$ with $\lambda \neq \mu.$
\item[(ii)]$
\displaystyle \mathcal{L}_2(g^- *_2 (S\otimes S))(\lambda,\mu) = \frac{1}{\lambda +\mu }[\hat g(\lambda) +\hat g(\mu)]\hat S(\lambda) \hat S(\mu)$ for $\Re\,\lambda,\,\Re\,\mu>\omega$ with $ \Re(\lambda+\mu)>0.$
\end{enumerate}
\end{proposition}
\bgproof
It is sufficient to apply Theorem \ref{convolution2} (iv), (i) (or (ii)) and \eqref{eq2.1}.
\edproof

The proof of the next corollary is a straightforward consequence of Theorem \ref{convolution2} (iv), Theorem \ref{cor3.1} and \eqref{eq2.1}.

\begin{corollary}\label{cor3.2}
Let $c\in L^1_{loc}(\R_+)$ be  an { absolutely continuous on $(0,\infty)$} and Laplace transformable function, and let $S:\R^+\to \mathcal{B}(X)$ be a locally integrable and strongly continuous and Laplace transformable operator valued function. \begin{enumerate}

\item[(i)]If $(c')^+:\R_+^2\to \CC$ is 2-Laplace transformable, then $$\displaystyle \mathcal{L}_2((c')^+ *_2 (S\otimes S))(\lambda,\mu) = \frac{1}{\mu - \lambda }[\lambda \hat c(\lambda) - \mu \hat c(\mu)]\hat S(\lambda) \hat S(\mu), \quad \Re\lambda,\Re\mu>\omega,\ \lambda \neq \mu.$$
\item[(ii)] If $(c')^-:\R_+^2\to \CC$ is 2-Laplace transformable and $c(0^+)$ exists then $$
\displaystyle \mathcal{L}_2((c')^- *_2 (S\otimes S))(\lambda,\mu) = \frac{1}{\lambda + \mu }[\lambda \hat c(\lambda) + \mu \hat c(\mu)]\hat S(\lambda) \hat S(\mu) - \frac{2c(0^+)}{\lambda + \mu}\hat S(\lambda) \hat S(\mu), $$ for $\Re \lambda,\Re\mu>\omega$ with $ \Re(\lambda+\mu)>0.$
\end{enumerate}
\end{corollary}

\begin{example}\label{ex3.3} {\rm
Let $c=g_{1-\alpha}$ for $0<\alpha<1$ and $\displaystyle \hat c(\lambda)= \frac{1}{\lambda^{1-\alpha}}$ for $\Re\lambda>0$. Then $c'= g_{-\alpha}$ and we obtain  the following identity,
\begin{equation}
\mathcal{L}_2((c')^+ *_2 (S\otimes S))(\lambda,\mu) = \frac{1}{\mu - \lambda }[\lambda^{\alpha} - \mu^{\alpha} ]\hat S(\lambda) \hat S(\mu),
\end{equation}
by Corollary \ref{cor3.2} (i). In particular, we recover \cite[Formula (2.8)]{Peng-Li-12}.}
\end{example}

\section{Local $(a,k)$-regularized resolvent families}
\label{local}
\setcounter{theorem}{0} \setcounter{equation}{0}

In this section we prove extension theorems for local $(a,k)$-regularized resolvent families. In the following we suppose that the function $k$ satisfies that $k(t)\neq 0$ for all $t\in (0,\sigma),$ where $\sigma$ is a sufficiently small positive number. We begin by { recalling} the following definition.

\begin{definition}\label{ak} Let $0<\tau\leq \infty,$ $a,k\in L^1_{loc}([0,\tau))$ with $k\in\mathcal{C}(0,\tau)$ {  that  $k(t)\neq 0$ for all $t\in (0,\sigma)$ ($\sigma$ small) } and $A$ be a closed operator. A strongly continuous operator family $\{S(t)\}_{t\in(0,\tau)}\subset \mathcal{B}(X)$  is a local (resp. global in case $\tau=\infty$) $(a,k)$-regularized resolvent family generated by $A$ if the following conditions are satisfied: \begin{itemize}

\item[(i)]$\displaystyle\lim_{t\to 0^+}\frac{S(t)x}{k(t)}=x$ for all $x\in X;$ \\

\item[(ii)]$S(t)A\subset AS(t),\quad {t\in (0,\tau)} $; \\

\item[(iii)]$(a*S)(t)x\in D(A)$ for $t\in(0,\tau)$ and $x\in X,$ and the following Volterra equation holds \begin{equation}\label{regularize}
A(a*S)(t)x=S(t)x-k(t)x,  \qquad x\in X,\,t\in(0,\tau).
\end{equation}
\end{itemize}
\end{definition}

\begin{remark}{\rm The reason why we do not consider directly the value of $S(\cdot)$ at the origin is that  $k$ could have a singularity at the origin; for example, $k(t)=g_{\beta+1}(t)$ has a singularity at $0$ if $-1<\beta<0.$ }
\end{remark}

\noindent In the rest of the paper we will assume that the functions $a,k$ are positive.

\noindent We note that loss of regularity arises because we treat with evolution equations corresponding to regularization of certain base equation. The typical example is the local $\alpha$-times integrated semigroups, e.g. the evolution equation:
$$
u'(t)= Au(t) + g_{\alpha}(t)x
$$
where  $\alpha > 0.$ In this case, the base equation, where no loss of regularity happens, is the Cauchy problem:
$$
u'(t)=Au(t),
$$
where $A$ is the generator of a $C_0$-semigroup, i.e.  it is known that a local $C_0$-semigroup can be extended without loss of regularity. In terms of $(a,k)$-regularized resolvent families, it means that if $A$ is the generator of a local $(1,1)$-regularized resolvent family on $[0,\tau)$ , then $A$ is also the generator of a $(1,1)$-regularized resolvent family on $[0,2\tau)$ and so on. However, this property is not { longer } true in the general case of $(a,k)$-regularized resolvent families, where loss of regularity is present. This phenomena has been observed for the case of $k$-convoluted semigroups \cite{Ke-Mi-La14} (in particular for $\alpha$-times integrated semigroups in \cite{Arendt-Elmen-Ke, Miana}) and $k$-convoluted cosine families \cite{Mi-Po14}.

In the next theorem we show the most general result about extension of $(a,k)$-regularized resolvent families. It shows that we can extend a local $(a,k)$-regularized resolvent family defined in $(0,\tau)$ to get a $(a,(k*a)^{*n}*a)$-regularized resolvent family in { $(0,(n+1)\tau).$}

\begin{theorem}\label{local1} Let $n\in \mathbb{N}$, $0<\tau\leq \infty$, $a,k\in L^1_{loc}([0,(n+1)\tau))$ with $k\in\mathcal{C}(0,(n+1)\tau),$ and $\{S_1(t)\}_{t\in(0,\tau)}$ be a local $(a,k)$-regularized resolvent family generated by $A.$ Then the family of operators $\{S_{n+1}(t)\}_{t\in (0,(n+1)T]}$ defined { recursively } by $$ S_{n+1}(t)x:=(k*a*S_n)(t)x, \qquad x\in X,\,t\in (0,nT],\text{ and}$$
$$
\begin{array}{lcl}
S_{n+1}(t)x &:=& \left(a^+\ast_2(S_n\otimes S_1)\right)(nT, t-nT)x +\displaystyle \int_0^{nT}k(t-r)(a*S_n)(r)x\,dr \\ \\ &\,&\qquad+ \displaystyle \int_0^{t-nT}((k*a)^{*(n-1)}*k)(t-r)(a*S_1)(r)x\,dr,
\end{array}
$$
for $x\in X$ and $t\in (nT,(n+1)T],$ is a local $(a,(k*a)^{*n}*k )$-regularized resolvent family generated by $A$ for any $T<\tau.$ Then { $A$ generates } a local $(a,(k*a)^{*n}*k )$-regularized resolvent family $\{S_{n+1}(t)\}_{t\in (0,(n+1)\tau)}.$

\end{theorem}
\bgproof Note that the family $\{S_{n+1}(t)\}_{t\in (0,(n+1)T]}$ is strongly continuous: \begin{displaymath}\begin{array}{l}
\displaystyle\lim_{t\to (nT)^+}S_{n+1}(t)x=\displaystyle\lim_{t\to (nT)^+}\biggl(\left(a^+\ast_2(S_n\otimes S_1)\right)(nT, t-nT)x \\ \\
+\displaystyle \int_0^{nT}k(t-r)(a*S_n)(r)x\,dr+ \displaystyle \int_0^{t-nT}((k*a)^{*(n-1)}*k)(t-r)(a*S_1)(r)x\,dr\biggr).
\end{array}\end{displaymath}
The first summand tends to $0$ using Corollary \ref{cor2.3} (i) in the vectorial case, the second one tends to $(k*a*S_n)(nT)x$, and the last term goes to $0.$ Furthermore, for $\varepsilon>0$ there exists $t>0$ sufficiently small such that $$\lVert \frac{S_n(s)x}{(k*a)^{*(n-1)}*k(s)}-x\rVert <\varepsilon, \text{ for }0<s<t.$$ Then  \small{\begin{displaymath}\begin{array}{l}
\displaystyle\lVert \frac{S_{n+1}(t)x}{(k*a)^{*n}*k(t)}-x\rVert=\lVert \frac{S_{n+1}(t)x-(k*a)^{*n}*k(t)x}{(k*a)^{*n}*k(t)}\rVert \\ \\
\displaystyle\leq \frac{1}{(k*a)^{*n}*k(t)}\int_0^{t}(k*a)(t-s)\lVert S_{n}(s)x-(k*a)^{*(n-1)}*k(s)x\rVert\,ds \\ \\
\displaystyle\leq \frac{1}{(k*a)^{*n}*k(t)}\int_0^{t}(k*a)(t-s)(k*a)^{*(n-1)}*k(s)\lVert \frac{S_{n}(s)x}{(k*a)^{*(n-1)}*k(s)}-x\rVert\,ds \leq \varepsilon.
\end{array}
\end{displaymath}}So, $\displaystyle\lim_{t\to 0^+}\frac{S_{n+1}(t)x}{(k*a)^{*n}*k(t)}=x$ for all $x\in X.$ Note  that $\{S_{n+1}(t)\}_{t\in (0,nT]}$  is a local $(a,(k*a)^{*n}*k )$-regularized resolvent family generated by $A$, see \cite[Remark 2.4 (4)]{Liz}. Now let $t\in (nT,(n+1)T]$ and $x\in X.$ It is clear that $S_{n+1}(t)A\subset AS_{n+1}(t).$ We show that $(a*S_{n+1})(t)x\in D(A).$ Note $$(a*S_{n+1})(t)x=\int_0^{nT}a(t-s)S_{n+1}(s)x\,ds+\int_{nT}^ta(t-s)S_{n+1}(s)x\,ds.$$ On the one hand, $$\int_0^{nT}a(t-s)S_{n+1}(s)x\,ds=\int_0^{nT}a(t-s)(k*a*S_{n})(s)x\,ds\in D(A),$$ since  $(a*S_n)(s)x\in D(A).$ On the other hand,

\small{\begin{eqnarray}
&&\nonumber\int_{nT}^{t}a(t-s)S_{n+1}(s)x\,ds \\
&&\nonumber= \displaystyle \int_{nT}^{t}a(t-s)\biggl(
\displaystyle \int_0^{s-nT}\int_0^{nT}a(s-r_1-r_2)S_n(r_1)S_1(r_2)x\,dr_1\,dr_2\biggr)\,ds
\\
\label{eq5.2} &&+  \displaystyle \int_{nT}^{t}a(t-s)\displaystyle \int_0^{nT}k(s-r)(a*S_n)(r)x\,dr\,ds \\\nonumber \\\nonumber &&+ \displaystyle \int_{nT}^{t}a(t-s)\displaystyle \int_0^{s-nT}((k*a)^{*(n-1)}*k)(s-r)(a*S_1)(r)x\,dr
\,ds.
\end{eqnarray}}

Note that $(a*S_n)(r)x, (a*S_1)(r)x\in D(A).$ Finally \small{\begin{eqnarray*}
&&\displaystyle \int_{nT}^{t}a(t-s)
\displaystyle \int_0^{s-nT}\int_0^{nT}a(s-r_1-r_2)S_n(r_1)S_1(r_2)x\,dr_1\,dr_2 \,ds \\
&&=\int_0^{nT}S_n(r_1)\int_{nT}^{t}a(t-s)\int_{nT}^s a(u-r_1)S_1(s-u)x\,du\,ds\,dr_1 \\
&&=\int_0^{nT}S_n(r_1)\int_{nT}^{t}a(u-r_1)\int_{u}^t a(t-s)S_1(s-u)x\,ds\,du\,dr_1 \\
&&=\int_0^{nT}S_n(r_1)\int_{nT}^{t}a(u-r_1)\int_0^{t-u}a(t-u-v)S_1(v)x\,dv\,du\,dr_1\in D(A)
\end{eqnarray*}} since $(a*S_1)(t-u)\in D(A).$
To finish the proof, we prove that for $t\in(nT,(n+1)T]$ and $x\in X$ the
equality \eqref{regularize} is satisfied. First observe that
$$A(a*S_{n+1})(t)x=A\int_0^{nT}a(t-s)(k*a*S_n)(s)x\,ds+A\int_{nT}^{t}a(t-s)S_{n+1}(s)x\,ds.$$
We apply the operator $A$ to the first summand of \eqref{eq5.2}, and
we obtain that \small{\begin{eqnarray*}
&A&\int_0^{nT}S_n(r_1)\int_{nT}^{t}a(t-s)\int_0^{s-nT}a(s-r_1-r_2)S_1(r_2)x\,dr_2\,ds\,dr_1 \\ &=&A\int_0^{nT}S_n(r_1)\int_{nT}^{t}a(t-s)\int_{nT}^s a(u-r_1)S_1(s-u)x\,du\,ds\,dr_1 \\
&=&A\int_0^{nT}S_n(r_1)\int_{nT}^{t}a(u-r_1)\int_{u}^t a(t-s)S_1(s-u)x\,ds\,du\,dr_1 \\
&=&A\int_0^{nT}S_n(r_1)\int_{nT}^{t}a(u-r_1)\int_0^{t-u}a(t-u-v)S_1(v)x\,dv\,du\,dr_1 \\
&=&\int_0^{nT}S_n(r_1)\int_{nT}^{t}a(u-r_1)\biggl(  S_1(t-u)-k(t-u) \biggr)x\,du\,dr_1 \\
&=&\int_0^{nT}S_n(r_1)\int_{0}^{t-nT}a(t-r_1-r_2)\biggl(
S_1(r_2)-k(r_2) \biggr)x\,dr_2\,dr_1
\end{eqnarray*}}
In the second summand of \eqref{eq5.2} we write
\small{$$
\begin{array}{l}
\displaystyle \int_0^{nT}(a*S_n)(r)x\int_{nT}^{t}a(t-s)k(s-r)\,ds\,dr \\
=\displaystyle \int_0^{nT}(a*S_n)(r)x\biggl(
\int_0^{t-r}a(t-r-u)k(u)\,du- \displaystyle \int_0^{nT-r}a(t-r-u)k(u)\,du\biggr)\,dr.
\end{array}
$$}
We apply the operator $A$ to each of the above terms to get
\small{\begin{eqnarray*}
&A& \int_0^{nT}(a*S_n)(r)x\int_0^{t-r}a(t-r-u)k(u)\,du\,dr \\ &=&A\int_0^{t-nT}k(u)\int_0^{nT}a(t-u-r)(a*S_n)(r)x\,dr\,du \\
&&+A\int_{t-nT}^t k(u)\int_0^{t-u}a(t-u-r)(a*S_n)(r)x\,dr\,du \\
&=&\int_0^{t-nT}k(u)\int_0^{nT}a(t-u-r)\biggl(S_n(r)x-((k*a)^{*(n-1)}*k)(r)x\biggr)\,dr\,du \\
&&+\int_{t-nT}^t k(u)\biggl( (a*S_n)(t-u)x-(k*a)^{*n}(t-u)x\biggr)\,du \\
&=&\int_0^{t-nT}k(u)\int_0^{nT}a(t-u-r)\biggl(S_n(r)x-((k*a)^{*(n-1)}*k)(r)x\biggr)\,dr\,du \\
&&+\int_{0}^{nT} k(t-r)\biggl(
(a*S_n)(r)x-(k*a)^{*n}(r)x\biggr)\,dr,
\end{eqnarray*}}
and \small{$$\begin{array}{l}
\displaystyle A\int_0^{nT}(a*S_n)(r)x\int_0^{nT-r}a(t-r-u)k(u)\,du\,dr\\
\displaystyle=A\int_0^{nT}k(u)\int_0^{nT-u}a(t-u-r)(a*S_n)(r)x\,dr\,du.
\end{array}$$}
In the third summand of \eqref{eq5.2} we write \small{\begin{eqnarray*} && \int_{0}^{t-nT}(a*S_1)(r)x\int_{r+nT}^{t}a(t-s)((k*a)^{*(n-1)}*k)(s-r)\,ds\,dr \\ &=&\int_{0}^{t-nT}(a*S_1)(r)x\biggl( \int_r^{t}a(t-s)((k*a)^{*(n-1)}*k)(s-r)\,ds \\
&&-\int_r^{r+nT}a(t-s)((k*a)^{*(n-1)}*k)(s-r)\,ds\biggr)\,dr.\end{eqnarray*}}
We apply the operator $A$ to each of the above terms to obtain
\small{\begin{eqnarray*}
&A& \int_{0}^{t-nT}(a*S_1)(r)x\int_r^{t}a(t-s)((k*a)^{*(n-1)}*k)(s-r)\,ds\,dr \\ &=&A\int_{0}^{t-nT}(a*S_1)(r)x\int_0^{t-r}a(t-r-u)((k*a)^{*(n-1)}*k)(u)\,du\,dr \\
&=&A\int_{0}^{nT}((k*a)^{*(n-1)}*k)(u)\int_0^{t-nT}a(t-u-r)(a*S_1)(r)x\,dr\,du \\
&&+A\int_{nT}^{t}((k*a)^{*(n-1)}*k)(u)\int_0^{t-u}a(t-u-r)(a*S_1)(r)x\,dr\,du 
\end{eqnarray*}}
\small{\begin{eqnarray*}
&=&A\int_{0}^{nT}((k*a)^{*(n-1)}*k)(u)\int_0^{t-nT}a(t-u-r)(a*S_1)(r)x\,dr\,du \\
&&+\int_{nT}^{t}((k*a)^{*(n-1)}*k(u))\biggl( (a*S_1)(t-u)x-(a*k)(t-u)x \biggr)\,du\\
&=&A\int_{0}^{nT}((k*a)^{*(n-1)}*k)(u)\int_0^{t-nT}a(t-u-r)(a*S_1)(r)x\,dr\,du \\
&&+\int_{0}^{t-nT}((k*a)^{*(n-1)}*k)(t-s)\biggl( (a*S_1)(s)x-(a*k)(s)x \biggr)\,ds\\
\end{eqnarray*}}
and
\small{\begin{eqnarray*}
&A&\int_0^{t-nT}(a*S_1)(r)x\int_r^{r+nT}a(t-s)((k*a)^{*(n-1)}*k)(s-r)\,ds\,dr\\ &=& A\int_0^{nT}((k*a)^{*(n-1)}*k)(u)\int_0^{t-nT}a(t-u-r)(a*S_1)(r)x\,dr\,du.
\end{eqnarray*}}
Then we have that
\small{\begin{eqnarray*}
&A&\int_{0}^{t-nT}(a*S_1)(r)x\int_{r+nT}^{t}a(t-s)((k*a)^{*(n-1)}*k)(s-r)\,ds\,dr \\&=& \int_0^{t-nT}((k*a)^{*(n-1)}*k)(t-s)\biggl((a*S_1)(s)x-(a*k)(s)x\biggr)\,ds.
\end{eqnarray*}}
Finally note that \small{\begin{eqnarray*}
&&A \int_0^{nT}a(t-s)(k*a*S_n)(s)x\,ds \\
&&=A\int_0^{nT}(a*S_n)(r)x\int_r^{nT}a(t-s)k(s-r)\,ds\,dr \\
&&=A\int_0^{nT}(a*S_n)(r)x\int_0^{nT-r}a(t-r-u)k(u)\,du\,dr \\
&&=A\int_0^{nT}k(u)\int_0^{nT-u}a(t-u-r)(a*S_n)(r)x\,dr\,du.
\end{eqnarray*}}
We join together all summands to conclude that
\small{$$\begin{array}{l}
\displaystyle A (a*S_{n+1})(t)x=\int_0^{t-nT}\int_0^{nT}a(t-r_1-r_2)S_n(r_1)S_1(r_2)x\,dr_1\,dr_2 \\  \\
\displaystyle + \int_0^{nT}k(t-r)(a*S_n)(r)x\,dr + \int_0^{t-nT}((k*a)^{*(n-1)}*k)(t-r)(a*S_1)(r)x\,dr \\ \\
\displaystyle  - \int_0^{t-nT}k(u)\int_0^{nT}a(t-u-r)((k*a)^{*(n-1)}*k)(r)x\,dr\,du\\ \\
\displaystyle - \int_0^{nT}k(t-r)(k*a)^{*n}(r)x\,dr- \int_0^{t-nT}((k*a)^{*(n-1)}*k)(t-r)(a*k)(r)x\,dr \\ \\
\displaystyle = S_{n+1}(t)x-(k*a)^{*n}(t)x,
\end{array}$$}
where we have used Lemma \ref{le5.1}.
\edproof

The expression of $\{S_{n+1}(t)\}_{t\in (0,(n+1)T]}$ is not unique; in the following theorem we show $\{S_{n+1}(t)\}_{t\in (0,(n+1)T]}$ in terms of $\{S_{j}(t)\}_{t\in (0,jT]}$ for all $1\leq j\leq n.$ The proof is similar to the proof of Theorem \ref{local1} and therefore we omit it.

\begin{theorem} Let $n\in \mathbb{N}$, $0<\tau\leq \infty$, $a, k\in L^1_{loc}([0,(n+1)\tau))$ with $k\in\mathcal{C}(0,(n+1)\tau),$ and $\{S_1(t)\}_{t\in (0,\tau)}$ be a local $(a,k)$-regularized resolvent family generated by $A.$ Then the family of operators $\{S_{n+1}(t)\}_{t\in (0,(n+1)T]}$ defined in Theorem \ref{local1} { satisfies }that \small{$$ S_{n+1}(t)x:=\biggl((k*a)^{*(n+1-j)}*S_j\biggr)(t)x, \qquad x\in X,\,t\in (0,jT],\text{ and}$$}
\small{$$\begin{array}{l}
S_{n+1}(t)x := \displaystyle  \left(a^+\ast_2 (S_j\otimes S_{n+1-j})\right)(jT, t-jT)x \\ \\
+ \displaystyle \int_0^{jT}((k*a)^{n-j}*k)(t-r)(a*S_j)(r)x\,dr \\ \\
+ \displaystyle \int_0^{t-jT}((k*a)^{*(j-1)}*k)(t-r)(a*S_{n+1-j})(r)x\,dr,
\end{array}$$}
for $x\in X,$ $1\leq j\leq n$ { and $t\in (jT,(n+1)T]$} for any $T<\tau.$

\end{theorem}
The following result is related to \cite[Theorem 4.4]{Ke-Mi-La14} and \cite[Theorem 3.3]{Mi-Po14}. However, note that both results are not included in this corollary.

\begin{corollary} Let $n\in \mathbb{N}$, $0<\tau\leq \infty$ and $\{S_1(t)\}_{t\in (0,\tau)}$ be a local $(g_{\alpha},g_{\beta+1})$-regularized resolvent family generated by $A,$ with $\alpha>0$ and $\beta>-1.$ Then the family of operators $\{S_{n+1}(t)\}_{t\in (0,(n+1)T]}$ defined by $$ S_{n+1}(t)x:=(g_{\beta+\alpha+1}*S_n)(t)x, \qquad x\in X, $$ for $t\in (0,nT]$ and
$$
\begin{array}{lcl}
S_{n+1}(t)x &:=& \displaystyle \left(g_\alpha^+\ast_2 (S_n\otimes S_1)\right)(nT,t-nT)x + \displaystyle \int_0^{nT}g_{\beta+1}(t-r)(g_{\alpha}*S_n)(r)x\,dr \\ \\
&\quad& + \displaystyle \int_0^{t-nT}g_{n(\beta+1)+\alpha(n-1)}(t-r)(g_{\alpha}*S_1)(r)x\,dr,
\end{array}
$$
for $x\in X$ { and $t\in (nT,(n+1)T]$} is a local $(g_\alpha,g_{(n+1)(\beta+1)+n\alpha} )$-regularized resolvent family generated by $A$ for any $T<\tau.$ Then $A$ generates a local $(g_\alpha,g_{(n+1)(\beta+1)+n\alpha} )$-regularized resolvent family $\{S_{n+1}(t)\}_{t\in (0,(n+1)\tau)}.$
\end{corollary}

However, if we restrict for example to the $\alpha$-times integrated semigroup case, the above extension is not the sharpest extension. Then for certain cases of the functions $a$ and $k$ there exist sharper extensions from the point of view of the regularized Cauchy problems. The following theorem gives us this sharp extension for a class of $(a,k)$-regularized resolvent families. Although the idea of the proof is similar to the proof of Theorem \ref{local1}, we have included it to make easier the reading because we use additional methods.

\begin{theorem}\label{local3} Let $n\in \mathbb{N}$, $0<\tau\leq \infty$, $a, k\in L^1_{loc}(\RR_+)$ with $k\in\mathcal{C}(0,\infty),$ Laplace transformable functions such that there exist $b,c\in L^1_{loc}(\R_+)$ Laplace transformable satisfying that $c$ is { absolutely continuous on $(0,\infty)$,} $(c')^+$ is 2-Laplace transformable,
$$(a*b)(t)=k(t),\ (a*c)(t)=1, \qquad t>0,
$$
and $\{S_1(t)\}_{t\in (0,\tau)}$ be a local $(a,k)$-regularized resolvent family generated by $A.$ Then the family of operators $\{S_{n+1}(t)\}_{t\in (0,(n+1)T]}$ defined by $$ S_{n+1}(t)x:=(b*S_n)(t)x, \qquad x\in X,$$ for $t\in (0,nT]$ and
\begin{eqnarray*}
S_{n+1}(t)x &:=& \displaystyle\int_0^{nT}b(t-r)S_n(r)x\,dr + \displaystyle \int_0^{t-nT}b^{*n}(t-r)S_1(r)x\,dr \\ \\
&&\displaystyle -\left((c')^+\ast_{2}(S_n\otimes S_1)\right)(nT,t-nT)x
\end{eqnarray*}
for $x\in X$ { and $t\in (nT,(n+1)T],$ } is a local $(a,b^{*n}*k )$-regularized resolvent family generated by $A$ for any $T<\tau.$ Then $A$ generates a local $(a,b^{*n}*k )$-regularized resolvent family $\{S_{n+1}(t)\}_{t\in (0,(n+1)\tau)}.$

\end{theorem}
\bgproof   Similarly to the proof of Theorem \ref{local1}, $\displaystyle\lim_{t\to 0^+}\frac{S_{n+1}(t)x}{(b^{*n}*k)(t)}=x$ for $x\in X.$ Note that $\{S_{n+1}(t)\}_{t\in (0,nT]}$  is a local $(a,b^{*n}*k )$-regularized resolvent family generated by $A$, see again \cite[Remark 2.4 (4)]{Liz}. Now let $t\in (nT,(n+1)T]$ and $x\in X.$ It is clear that $S_{n+1}(t)A\subset AS_{n+1}(t),$ and following the proof of  Theorem \ref{local1} it is easy to see that $(a*S_{n+1})(t)x\in D(A)$.

Now we prove that for $t\in(nT,(n+1)T]$ and $x\in X$ the
equality \eqref{regularize} is satisfied. First observe that
\small{$$A(a*S_{n+1})(t)x=A\int_0^{nT}a(t-s)(b*S_n)(s)x\,ds+A\int_{nT}^{t}a(t-s)S_{n+1}(s)x\,ds.$$}
Note that
\small{\begin{equation}
\begin{array}{l}
\displaystyle \int_{nT}^{t}a(t-s)S_{n+1}(s)x\,ds = \displaystyle\int_{nT}^{t}a(t-s)\biggl(\displaystyle \int_0^{nT}b(s-r)S_n(r)x\,dr \\ \\
 \label{eq7.1} + \displaystyle \int_0^{s-nT}b^{*n}(s-r)S_1(r)x\,dr \biggr)\,ds\\ \\
\displaystyle -\int_{nT}^{t}a(t-s)\int_0^{s-nT}\int_0^{nT}c'(s-r_1-r_2)S_n(r_1)S_1(r_2)x\,dr_1\,dr_2\,ds.
\end{array}
\end{equation}}
We apply the operator $A$ to the third summand of \eqref{eq7.1}, and
we obtain that \small{\begin{eqnarray*}
&-A&\int_0^{nT}S_n(r_1)\int_{nT}^{t}a(t-s)\int_0^{s-nT}c'(s-r_1-r_2)S_1(r_2)x\,dr_2\,ds\,dr_1 \\ \\
&=&-A\int_0^{nT}S_n(r_1)\int_{nT}^{t}a(t-s)\int_{nT}^s c'(u-r_1)S_1(s-u)x\,du\,ds\,dr_1 \\ \\
&=&-A\int_0^{nT}S_n(r_1)\int_{nT}^{t}c'(u-r_1)\int_{u}^t a(t-s)S_1(s-u)x\,ds\,du\,dr_1 
\end{eqnarray*}}
\small{\begin{eqnarray*}
&=&-A\int_0^{nT}S_n(r_1)\int_{nT}^{t}c'(u-r_1)\int_0^{t-u}a(t-u-v)S_1(v)x\,dv\,du\,dr_1 \\ \\
&=&-\int_0^{nT}S_n(r_1)\int_{nT}^{t}c'(u-r_1)\biggl(  S_1(t-u)-k(t-u) \biggr)x\,du\,dr_1 \\ \\
&=&-\int_0^{nT}S_n(r_1)\int_{0}^{t-nT}c'(t-r_1-r_2)\biggl(
S_1(r_2)-k(r_2) \biggr)x\,dr_2\,dr_1.
\end{eqnarray*}}
In the first summand of \eqref{eq7.1} we write
\small{$$
\begin{array}{l}
\displaystyle \int_0^{nT}S_n(r)x\int_{nT}^{t}a(t-s)b(s-r)\,ds\,dr =\displaystyle \int_0^{nT}S_n(r)x
\int_0^{t-r}a(t-r-u)b(u)\,du\,dr\\ - \displaystyle \int_0^{nT}S_n(r)x \int_0^{nT-r}a(t-r-u)b(u)\,du\,dr.
\end{array}
$$}
We apply the operator $A$ to each of the above terms to get
\small{\begin{eqnarray*}
&A& \int_0^{nT}S_n(r)x\int_0^{t-r}a(t-r-u)b(u)\,du\,dr \\ &=&A\int_0^{t-nT}b(u)\int_0^{nT}a(t-u-r)S_n(r)x\,dr\,du \\ \\
&&+A\int_{t-nT}^t b(u)\int_0^{t-u}a(t-u-r)S_n(r)x\,dr\,du \\ \\
&=&A\int_0^{t-nT}b(u)\int_0^{nT}a(t-u-r)S_n(r)x\,dr\,du \\ \\
&&+\int_{t-nT}^t b(u)\biggl( S_n(t-u)x-(b^{*(n-1)}*k)(t-u)x\biggr)\,du \\ \\
&=&A\int_0^{t-nT}b(u)\int_0^{nT}a(t-u-r)S_n(r)x\,dr\,du \\ \\
&&+\int_{0}^{nT} b(t-r)\biggl(S_n(r)x-(b^{*(n-1)}*k)(r)x\biggr)\,dr,
\end{eqnarray*}}
and 
\small{$$\int_0^{nT}S_n(r)x\int_0^{nT-r}a(t-r-u)b(u)\,du\,dr 
=\int_0^{nT}b(u)\int_0^{nT-u}a(t-u-r)S_n(r)x\,dr\,du.$$}
In the second summand of \eqref{eq7.1} we write \small{\begin{eqnarray*} && \int_{0}^{t-nT}S_1(r)x\int_{r+nT}^{t}a(t-s)b^{*n}(s-r)\,ds\,dr \\ \\
 &=&\int_{0}^{t-nT}S_1(r)x\biggl( \int_r^{t}a(t-s)b^{*n}(s-r)\,ds-\int_r^{r+nT}a(t-s)b^{*n}(s-r)\,ds\biggr)\,dr.\end{eqnarray*}}
We apply the operator $A$ to each of the above terms to obtain
\small{\begin{eqnarray*}
&A& \int_{0}^{t-nT}S_1(r)x\int_r^{t}a(t-s)b^{*n}(s-r)\,ds\,dr \\ \\
&=& A\int_{0}^{t-nT}S_1(r)x\int_0^{t-r}a(t-r-u)b^{*n}(u)\,du\,dr \\ \\
&=&A\int_{0}^{nT}b^{*n}(u)\int_0^{t-nT}a(t-u-r)S_1(r)x\,dr\,du \\ \\
&&+ A\int_{nT}^{t}b^{*n}(u)\int_0^{t-u}a(t-u-r)S_1(r)x\,dr\,du \\ \\
&=&A\int_{0}^{nT}b^{*n}(u)\int_0^{t-nT}a(t-u-r)S_1(r)x\,dr\,du \\ \\
&&+ \int_{nT}^{t}b^{*n}(u)\biggl(S_1(t-u)x-k(t-u)x \biggr)\,du\\ \\
&=&A\int_{0}^{nT}b^{*n}(u)\int_0^{t-nT}a(t-u-r)S_1(r)x\,dr\,du \\ \\
&&+ \int_{0}^{t-nT}b^{*n}(t-s)\biggl( S_1(s)x-k(s)x \biggr)\,ds
\end{eqnarray*}}
and
$$
\int_0^{t-nT}S_1(r)x\int_r^{r+nT}a(t-s)b^{*n}(s-r)\,dsdr 
 =\int_0^{nT}b^{*n}(u)\int_0^{t-nT}a(t-u-r)S_1(r)x\,drdu.
 $$
Then we have that
\small{$$
A\int_{0}^{t-nT}S_1(r)x\int_{r+nT}^{t}a(t-s)b^{*n}(s-r)\,ds\,dr= \int_0^{t-nT}b^{*n}(t-s)\biggl(S_1(s)x-k(s)x\biggr)\,ds.
$$}
Furthermore note that \small{\begin{eqnarray*}
&\,&A \int_0^{nT}a(t-s)(b*S_n)(s)x\,ds=A\int_0^{nT}S_n(r)x\int_r^{nT}a(t-s)b(s-r)\,ds\,dr \\ \\
&\quad&=A\int_0^{nT}S_n(r)x\int_0^{nT-r}a(t-r-u)b(u)\,du\,dr \\ \\
 &\quad&=A\int_0^{nT}b(u)\int_0^{nT-u}a(t-u-r)S_n(r)x\,dr\,du.
\end{eqnarray*}}
We join together all summands { to conclude} that
\small{\begin{eqnarray*}
A (a*S_{n+1})(t)x&=&S_{n+1}(t)x+ A\int_0^{t-nT}b(u)\int_0^{nT}a(t-u-r)S_n(r)x\,dr\,du \\
&-& \int_0^{nT}b(t-r)(b^{*(n-1)}*k)(r)x\,dr- \int_0^{t-nT}b^{*n}(t-r)k(r)x\,dr \\
&+&\int_0^{nT}S_n(r_1)x\int_{0}^{t-nT}c'(t-r_1-r_2)k(r_2)\,dr_2\,dr_1.
\end{eqnarray*}}
Now we use induction. As $\{S_n(t)\}_{t\in (0,nT]}$ is a local $(a,b^{*(n-1)}*k)$-regularized resolvent family generated by $A,$ then $$S_n(r_1)x=A(a*S_n)(r_1)x+(b^{*(n-1)}*k)(r_1)x=A(a*S_n)(r_1)x+(b^{*n}*a)(r_1)x,$$ and so

\small{\begin{displaymath}
\begin{array}{lcl}
&&\displaystyle \int_0^{nT}S_n(r_1)x\int_{0}^{t-nT}c'(t-r_1-r_2)k(r_2)\,dr_2\,dr_1 \\ \\
&=&\displaystyle A\int_0^{nT}(a*S_n)(r_1)x\int_{0}^{t-nT}c'(t-r_1-r_2)(a*b)(r_2)\,dr_2\,dr_1 \\ \\
&\quad&+\displaystyle \int_0^{nT}(b^{*n}*a)(r_1)x\int_{0}^{t-nT}c'(t-r_1-r_2)(a*b)(r_2)\,dr_2\,dr_1.
\end{array}
\end{displaymath}}
On one hand, \small{\begin{displaymath}
\begin{array}{l}
\displaystyle A\int_0^{nT}(a*S_n)(r_1)x\int_{0}^{t-nT}c'(t-r_1-r_2)(a*b)(r_2)\,dr_2\,dr_1 \\ \\
=\displaystyle A\int_0^{nT}S_n(u)x\int_{0}^{t-nT}b(v)\int_u^{nT}\int_v^{t-nT}a(r_1-u)a(r_2-v)c'(t-r_1-r_2)\,dr_2\,dr_1\,dv\,du \\ \\
=\displaystyle -A\int_0^{nT}S_n(u)x\int_{0}^{t-nT}a(t-u-v)b(v)\,dv\,du
\end{array}
\end{displaymath}}
and on the other hand \small{\begin{displaymath}
\begin{array}{l}
\displaystyle \int_0^{nT}(b^{*n}*a)(r_1)x\int_{0}^{t-nT}c'(t-r_1-r_2)(a*b)(r_2)\,dr_2\,dr_1\\ \\
=\displaystyle\int_0^{nT}b^{*n}(u)x\int_{0}^{t-nT}b(v)\int_u^{nT}\int_v^{t-nT}a(r_1-u)a(r_2-v)c'(t-r_1-r_2)\,dr_2\,dr_1\,dv\,du \\ \\
=\displaystyle -\int_0^{nT}b^{*n}(u)x\int_{0}^{t-nT}a(t-u-v)b(v)\,dv\,du,
\end{array}
\end{displaymath}}
where we have applied Theorem \ref{le3.4} (i). Applying Lemma \ref{le5.1} we get that
\small{$$\begin{array}{l}
A (a*S_{n+1})(t)x=S_{n+1}(t)x-\displaystyle\int_0^{nT}b(t-r)(b^{*(n-1)}*k)(r)x\,dr \\ \\
- \displaystyle\int_0^{t-nT}b^{*n}(t-r)k(r)x\,dr - \displaystyle\int_0^{nT}b^{*n}(u)x\int_{0}^{t-nT}a(t-u-v)b(v)\,dv\,du \\ \\
=S_{n+1}(t)x-\displaystyle\int_0^{nT}b(t-r)(b^{*n}*a)(r)x\,dr \\ \\
- \displaystyle\int_0^{t-nT}b^{*n}(t-r)(a*b)(r)x\,dr - \displaystyle\int_0^{nT}b^{*n}(u)x\int_{0}^{t-nT}a(t-u-v)b(v)\,dv\,du \\ \\
= S_{n+1}(t)x-(b^{*(n+1)}*a)(t)x=S_{n+1}(t)x-(b^{*n}*k)(t)x.
\end{array}$$}

Finally we check that the family  $\{S_{n+1}(t)\}_{t\in (0,(n+1)T]}$ is strongly continuous. It is direct to check that $\{S_{n+1}(t)\}_{t\in (0, (n+1)T]}$ is uniformly bounded on  $[nT-\varepsilon,nT+\varepsilon]$ for all $0<\varepsilon<T,$ and strongly continuous on $(0, nT)\cup (nT, (n+1)T]$. Note that for $t\to (nT)^+$, we have that
\small{$$\begin{array}{l}
\displaystyle S_{n+1}(t)x-(b^{\ast n}\ast k)(t)x= A\left(\int_0^{nT}a(t-s)S_{n+1}(s)xds\right) \\ \\
\displaystyle +A\left(\int_{nT}^ta(t-s)S_{n+1}(s)xds\right)\to A\left(\int_0^{nT}a(nT-s)S_{n+1}(s)xds\right) \\ \\
\displaystyle = S_{n+1}(nT)x-(b^{\ast n}\ast k)(nT)x, \qquad x\in X,
\end{array}$$}
and we conclude that the family $\{S_{n+1}(t)\}_{t\in [0,(n+1)T]}$ is strongly continuous.\edproof

The following result extends   \cite[Theorem 2]{Miana}, because we obtain the sharp extension of $(g_{\alpha},g_{\beta+1})$-regularized resolvent families when $0<\alpha<1$ and $\beta-\alpha>-1,$ and when $\alpha\to 1^-$ we recover the  $\alpha$-times integrated semigroup case, considered in \cite{Miana}. More generally one could consider the case of $K$-convoluted resolvent families, i.e. $(g_{\alpha},(1*K))$-regularized resolvent families for $0<\alpha<1$ and compare it to the limit case when $\alpha\to 1^-,$ see \cite[Theorem 4.4]{Ke-Mi-La14}.

\begin{corollary}\label{Cor5.5} Let $n\in \mathbb{N}$, $0<\tau\leq \infty$ and $\{S_1(t)\}_{t\in (0,\tau)}$ be a local $(g_{\alpha},g_{\beta+1})$-regularized resolvent family generated by $A$ with $0<\alpha<1$ and $\beta-\alpha>-1.$ Then the family of operators $\{S_{n+1}(t)\}_{t\in (0,(n+1)T]}$ defined by $$ S_{n+1}(t)x:=(g_{\beta-\alpha+1}*S_n)(t)x, \qquad x\in X, $$ for $t\in (0,nT]$ and
$$
\begin{array}{lcl}
S_{n+1}(t)x &:=& \displaystyle\int_0^{nT}g_{\beta-\alpha+1}(t-r)S_n(r)x\,dr + \displaystyle \int_0^{t-nT}g_{n(\beta-\alpha+1)}(t-r)S_1(r)x\,dr \\ \\
&\,&-\left((g_{-\alpha})^+\ast \left(S_n\otimes S_1\right)\right)(nT,t-nT)x,
\end{array}
$$
for $x\in X$ { and $t\in (nT,(n+1)T],$ } is a local $(g_{\alpha},g_{n(\beta-\alpha+1)+\beta+1})$-regularized resolvent family generated by $A$ for any $T<\tau.$ Then $A$ generates a local $(g_{\alpha},g_{n(\beta-\alpha+1)+\beta+1})$-regularized resolvent family $\{S_{n+1}(t)\}_{t\in (0,(n+1)\tau)}.$
\end{corollary}

\section{Solutions of evolutionary problems without jumps of regularity}
\label{jump}
\setcounter{theorem}{0} \setcounter{equation}{0}

In this section, we identify  a wide class of evolution equations where no loss of regularity happens. It is interesting to note that it was not known until now if this property goes beyond the cases of the heat and wave equations, i.e., the semigroup and cosine cases.  We begin with the following result which is subordinated to the semigroup case in the sense that we cannot go beyond of $\alpha >1$ when we restrict to the particular case of $(g_{\alpha},g_{\alpha})$-regularized resolvent families. See the next corollary.

\begin{theorem}\label{th-local1} Let $n\in \mathbb{N}$, $0<\tau\leq \infty$, $a\in L^1_{loc}(\R_+)$ with $a\in\mathcal{C}(0,\infty),$ { be a }Laplace transformable function such that there exists $c\in L^1_{loc}(\R_+)$ Laplace transformable satisfying that $c$ is absolutely continuous on $(0,\infty)$, $(c')^+$ is 2-Laplace transformable and $(a*c)(t)=1$ for all $t>0,$ and $\{S_1(t)\}_{t\in (0,\tau)}$ be a local $(a,a)$-regularized resolvent family generated by $A.$ Then the family of operators $\{S_{n+1}(t)\}_{t\in (0,(n+1)T]}$ defined by
$$S_{n+1}(t)x:=S_n(t)x, \qquad x\in X, $$ for $t\in (0,nT]$ and
$$\begin{array}{c}
S_{n+1}(t)x := -((c')^+ *_2 (S_n\otimes S_1))(nT,t-nT)x, \qquad x\in X,\\ \\
\end{array}
$$
and $t\in (nT,(n+1)T]$ is a local $(a,a)$-regularized resolvent family generated by $A$ for any $T<\tau.$ { Then $A$ } generates a global $(a,a)$-regularized resolvent family $\{S(t)\}_{t\in (0,\infty)}.$
\end{theorem}
\bgproof
Note that $\displaystyle{\lim_{t\to 0^+}\frac{S_{n+1}(t)x}{a(t)}=x}$ for $x\in X$ and the family $\{S_{n+1}(t)\}_{t\in (0,(n+1)T]}$ is strongly continuous. The proof of this fact is similar to  Theorem \ref{local3}. Obviously, $\{S_{n+1}(t)\}_{t\in (0,nT]}$  is a local $(a,a)$-regularized resolvent family generated by $A$. Now let $t\in (nT,(n+1)T]$ and $x\in X.$ It is clear that $S_{n+1}(t)A\subset AS_{n+1}(t).$ We show that $(a*S_{n+1})(t)x\in D(A).$ Note $$(a*S_{n+1})(t)x=\int_0^{nT}a(t-s)S_{n}(s)x\,ds+\int_{nT}^t a(t-s)S_{n+1}(s)x\,ds, \qquad x\in X.$$ On  one hand, note that
$\displaystyle{\int_0^{nT}a(t-s)S_{n}(s)x\,ds\in D(A),}$ see \eqref{eq6.3} at the end of the proof. On the other hand,
\small{$$\begin{array}{l}
\displaystyle \int_{nT}^{t}a(t-s)S_{n+1}(s)x\,ds \\ \\
\displaystyle=-\int_{nT}^{t}a(t-s)\int_0^{s-nT}\int_0^{nT}c'(s-r_1-r_2)S_n(r_1)S_1(r_2)x\,dr_1\,dr_2\,ds \\ \\
\displaystyle =-\int_0^{nT}S_n(r_1)\int_{nT}^{t}a(t-s)\int_0^{s-nT}c'(s-r_1-r_2)S_1(r_2)x\,dr_2\,ds\,dr_1 \\ \\
\displaystyle =-\int_0^{nT}S_n(r_1)\int_{nT}^{t}a(t-s)\int_{nT}^s c'(u-r_1)S_1(s-u)x\,du\,ds\,dr_1 \\ \\
\displaystyle =-\int_0^{nT}S_n(r_1)\int_{nT}^{t}c'(u-r_1)\int_{u}^t a(t-s)S_1(s-u)x\,ds\,du\,dr_1 \\ \\
\displaystyle =-\int_0^{nT}S_n(r_1)\int_{nT}^{t}c'(u-r_1)\int_{0}^{t-u} a(t-u-v)S_1(v)x\,dv\,du\,dr_1 \in D(A)
\end{array}$$}
since $(a*S_1)(t-u)\in D(A).$ To finish the proof, we prove that for $t\in(nT,(n+1)T]$ and $x\in X$ the equality $$A(a*S_{n+1}(t))x=S_{n+1}(t)x-a(t)x,$$ is verified. First observe that
$$A(a*S_{n+1})(t)x=A\int_0^{nT}a(t-s)S_n(s)x\,ds+A\int_{nT}^{t}a(t-s)S_{n+1}(s)x\,ds.$$ Now, we develop the second term applying change of variables and Fubini's theorem:
\small{$$\begin{array}{l}
\displaystyle A\int_{nT}^{t}a(t-s)S_{n+1}(s)x\,ds \\ \\
\displaystyle = -A\int_{nT}^{t}a(t-s)\int_0^{s-nT}\int_0^{nT}c'(s-r_1-r_2)S_n(r_1)S_1(r_2)x\,dr_1\,dr_2\,ds \\ \\
\displaystyle = -A\int_0^{nT}S_n(r_1)\int_{nT}^{t}a(t-s)\int_0^{s-nT}c'(s-r_1-r_2)S_1(r_2)x\,dr_2\,ds\,dr_1 \\ \\
\displaystyle = -A\int_0^{nT}S_n(r_1)\int_{nT}^{t}a(t-s)\int_{nT}^s c'(u-r_1)S_1(s-u)x\,du\,ds\,dr_1 \\ \\
\displaystyle = -A\int_0^{nT}S_n(r_1)\int_{nT}^{t}c'(u-r_1)\int_{u}^t a(t-s)S_1(s-u)x\,ds\,du\,dr_1 \\ \\
\displaystyle = -A\int_0^{nT}S_n(r_1)\int_{nT}^{t}c'(u-r_1)\int_{0}^{t-u} a(t-u-v)S_1(v)x\,dv\,du\,dr_1 \\ \\
\displaystyle = -\int_0^{nT}S_n(r_1)\int_{nT}^{t}c'(u-r_1)(S_1(t-u)-a(t-u))x\,du\,dr_1 \\ \\
\displaystyle = -\int_0^{nT}S_n(r_1)\int_{0}^{t-nT}c'(t-r_1-r_2)(S_1(r_2)-a(r_2))x\,dr_2\,dr_1,
\end{array}$$}
where we have used that $\{S_1(t)\}_{t\in (0,T]}$ is a local $(a,a)$-regularized resolvent family generated by $A.$ Then \small{\begin{eqnarray*}
A(a*S_{n+1})(t)&=&A\int_0^{nT}a(t-s)S_n(s)x\,ds+ S_{n+1}(t)x \\ \\
&+&\int_0^{nT}S_n(r_1)x\int_{0}^{t-nT}c'(t-r_1-r_2)a(r_2)\,dr_2\,dr_1.
\end{eqnarray*}}
As $\{S_n(t)\}_{t\in (0,nT]}$ is a local $(a,a)$-regularized resolvent family generated by $A,$ then $$S_n(r_1)x=A(a*S_n)(r_1)x+a(r_1)x,$$ and \small{\begin{eqnarray*} &&\int_0^{nT}S_n(r_1)x\int_{0}^{t-nT}c'(t-r_1-r_2)a(r_2)\,dr_2\,dr_1 \\ \\
&&=\int_0^{nT}(A(a*S_n)(r_1)+a(r_1))x\int_{0}^{t-nT}c'(t-r_1-r_2)a(r_2)\,dr_2\,dr_1.
\end{eqnarray*}}
On the one hand, we obtain the following identity by change of variables and Fubini's theorem:
\small{\begin{equation}
\begin{array}{lcl}\label{eq6.3}
&&\displaystyle A\int_0^{nT}(a*S_n)(r_1)x\int_{0}^{t-nT}c'(t-r_1-r_2)a(r_2)\,dr_2\,dr_1 \\ \\
&&\displaystyle =A\int_0^{nT}(\int_0^{r_1}a(r_1-u)S_n(u)x\,du)\int_{0}^{t-nT}c'(t-r_1-r_2)a(r_2)\,dr_2\,dr_1 \\ \\
&&\displaystyle =A\int_0^{nT}S_n(u)x\int_u^{nT}a(r_1-u)\int_{0}^{t-nT}c'(t-r_1-r_2)a(r_2)\,dr_2\,dr_1\,du \\ \\
&&\displaystyle =A\int_0^{nT}S_n(u)x\int_0^{nT-u}\int_{0}^{t-nT}c'(t-u-v-r_2)a(v)a(r_2)\,dr_2\,dv\,du \\ \\
&&\displaystyle =A\int_0^{nT}S_n(u)x((c')^+ *_2 (a\otimes a))(nT-u,t-nT)\,du \\ \\
&&\displaystyle =-A\int_0^{nT}a(t-u)S_n(u)x\,du,
\end{array}
\end{equation}}
\noindent where we have used Theorem \ref{le3.4}. On the other hand, we use Theorem \ref{le3.4} again to get $$\begin{array}{l}
\displaystyle\int_0^{nT}a(r_1)x\int_{0}^{t-nT}c'(t-r_1-r_2)a(r_2)\,dr_2\,dr_1=-((c')^+ *_2 (a\otimes a ))(t-nT,nT)x\\ \\
=-a(t)x.\end{array}$$ We join all the terms and we obtain the result.
\edproof

The next result considers the special case of $(g_{\alpha},g_{\alpha})$-regularized resolvent families. Here we have to restrict to the range $0<\alpha<1$ according to the given hypothesis in the above Theorem. We observe that this condition is optimal in the following sense: When $\alpha=1$ we are treating with the parabolic case, i.e. the equation
\begin{equation*}
\left\{\begin{array}{ll}
u'(t)=Au(t)+x,&t \in [0,\tau), \quad  x \in D(A), \\
u(0)=0,&
\end{array} \right.
\end{equation*}
where $A$ is the generator of a $C_0$-semigroup, or, equivalently, a $(1,1)$-regularized resolvent family. We known that in this case no loss of regularity happens. Now,  for $0<\alpha<1$ we have to consider the fractional order differential equation:
\begin{equation}\label{eq6.1}
\left\{\begin{array}{ll}
_RD_t^{\alpha}u(t)= Au(t) + g_{\alpha}(t)x,&t \in (0,\tau), \quad  x \in D(A), \\
(g_{1-\alpha}*u)(0)=0,&
\end{array} \right.
\end{equation}
where $_RD_t$ denotes the fractional derivative in the Riemann Liouville sense, and $A$ is the generator of a $(g_{\alpha},g_{\alpha})$-regularized resolvent family { (see also \cite[Example 2.1.38]{Kostic2} and the paragraph preceding it).} The following corollary shows that again no loss of regularity happens for equation (\ref{eq6.1}). In passing, we conclude the remarkable fact that equation (\ref{eq6.1}) is at the basis of the process of regularization for $0<\alpha<1,$ where the solutions of the regularized problems correspond to the families of the Corollary \ref{Cor5.5}.

In the following picture we can see graphically the previous comments for $(g_{\alpha},g_{\beta+1})$-regularized resolvent families with $0<\alpha<1$. Note that the straight line formed by the points $(\alpha,\alpha-1)$ corresponds to $(g_{\alpha},g_{\alpha})$-regularized families, which is the basis of the process of regularization for $0<\alpha<1.$ For $\alpha=1,$ the point $(1,0)$ corresponds to a $C_0$-semigroup, and the points $(1,\beta)$ correspond to $\beta$-times integrated semigroups for $\beta>0.$

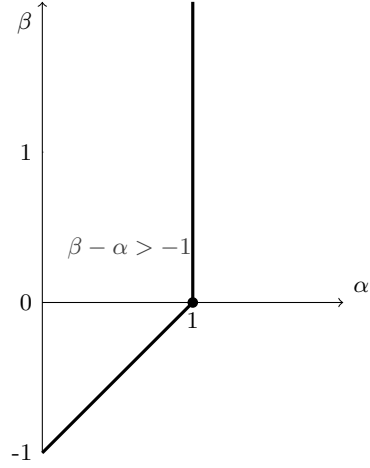
\begin{center}
\begin{tikzpicture}
   \begin{scope}[font=\scriptsize,scale=2]
     Axes:
     Are simply drawn using line with the `->` option to make them arrows:
    The main labels of the axes can be places using `node`s:
    \draw [->] (0,0) -- coordinate (x axis mid) (2,0) node [above right]  {$\alpha$};
   \draw [->] (0,-1) -- coordinate (y axis mid) (0,2) node [below left] {$\beta$};

    \foreach \x in {1}
     		\draw (\x,0) -- (\x,0)
			node[anchor=north] {\x};
    	\foreach \y in {-1,0,1}
     		\draw (0,\y) -- (0,\y)
     			node[anchor=east] {\y};

    \draw[very thick] (0,-1) coordinate (a_1) -- (1,0) coordinate (a_2);

    \draw[very thick] (1,0) coordinate (b_1) -- (1,2) coordinate (b_2);

    \coordinate (c) at (intersection of a_1--a_2 and b_1--b_2);

    \fill[black] (c) circle (1pt);

      (0,-1) -- (0,2) -- (1,2) -- (1,0) -- cycle;

    \node [below right,darkgray] at (0.1,0.5) {$\beta-\alpha>-1$};

    \end{scope}

\end{tikzpicture}
\end{center}

\begin{corollary}\label{co-local1} Let $n\in \mathbb{N}$, $0<\tau\leq \infty,$ $0<\alpha<1$ and $\{S_1(t)\}_{t\in (0,\tau)}$ be a local $(g_{\alpha},g_{\alpha})$-regularized resolvent family generated by $A.$ Then the family of operators $\{S_{n+1}(t)\}_{t\in (0,(n+1)T]}$ defined by

$$ S_{n+1}(t)x=S_n(t)x, \qquad x\in X, $$ for $t\in (0,nT]$ and

$$\begin{array}{c}
S_{n+1}(t)x = \displaystyle \frac{\alpha}{\Gamma(1-\alpha)}\int_0^{t-nT}\int_0^{nT}\frac{S_n(r_1)S_1(r_2)x}{(t-r_1-r_2)^{1+\alpha}}\,dr_1\,dr_2 \\ \\
\end{array}$$
for $x\in X$ { and $t\in (nT,(n+1)T]$} is a local $(g_{\alpha},g_{\alpha})$-regularized resolvent family generated by $A$ for any $T<\tau.$ { Then $A$} generates a global $(g_{\alpha},g_{\alpha})$-regularized resolvent family $\{S(t)\}_{t\in (0,\infty)}.$
\end{corollary}

Now, we consider a different class of $(a,k)$-regularized resolvent families such that we can solve the extension problem without loss of regularity.

\begin{theorem}\label{th-local2} Let $n\in \mathbb{N}$, $0<\tau\leq \infty$, $a\in L^1_{loc}(\R_+)$ with $a\in\mathcal{C}(0,\infty),$ Laplace transformable function such that there exists $c\in L^1_{loc}(\R_+)$ Laplace transformable satisfying that $c$ is absolutely continuous, differentiable a.e., $c(0^+)=0,$ $(c')^+$ and $(c')^-$ are 2-Laplace transformable, and $(a*c)(t)=1$ for all $t>0,$ and $\{S_1(t)\}_{t\in (0,\tau)}$ be a local $(a*1,a)$-regularized resolvent family generated by $A.$ Then the family of operators $\{S_{n+1}(t)\}_{t\in (0,(n+1)T]}$ defined by $$ S_{n+1}(t)x:=S_n(t)x, \qquad x\in X, $$ for $t\in (0,nT]$ and
\begin{eqnarray*}
S_{n+1}(t)x &:=&-S_n(2nT-t)x +((c')^- *_2 (S_n\otimes S_1))(nT,t-nT)x  \\ \\
&&- ((c')^+ *_2 (S_n\otimes S_1))(nT,t-nT)x
\end{eqnarray*}
for $x\in X$ { and $t\in (nT,(n+1)T]$} is a local $(a*1,a)$-regularized resolvent family generated by $A$ for any $T<\tau.$ Then $A$ generates a global $(a*1,a)$-regularized resolvent family $\{S(t)\}_{t\in (0,\infty)}.$
\end{theorem}

\bgproof Note that $\displaystyle{\lim_{t\to 0^+}\frac{S_{n+1}(t)x}{a(t)}=x}$ for $x\in X$ and the family   $\{S_{n+1}(t)\}_{t\in (0,(n+1)T]}$ is strongly continuous, see the proof in Theorem \ref{local3}; in particular, $\{S_{n+1}(t)\}_{t\in (0,nT]}$  is a local $(a*1,a)$-regularized resolvent family generated by $A$. Now let $t\in (nT,(n+1)T]$ and $x\in X.$ It is clear that $S_{n+1}(t)A\subset AS_{n+1}(t).$ Following the same ideas as in the proofs of the previous theorems, we conclude that $(a*1*S_{n+1})(t)x\in D(A).$

To finish the proof, it remains to prove that for $t\in(nT,(n+1)T]$ and $x\in X$ the equality $A(a*1*S_{n+1}(t))x=S_{n+1}(t)x-a(t)x,$ is satisfied. First observe that
$$A(a*1*S_{n+1})(t)x=A\int_0^{nT}(a*1)(t-s)S_n(s)x\,ds+A\int_{nT}^{t}(a*1)(t-s)S_{n+1}(s)x\,ds.$$ Note that
\small{$$\begin{array}{l}
\displaystyle \int_{nT}^{t}(a*1)(t-s)S_{n+1}(s)x\,ds = \int_{nT}^{t}(a*1)(t-s)\biggl( -S_n(2nT-t)x \\ \\
\displaystyle +\int_0^{s-nT}\int_0^{nT}(c')^-(nT-r_1,s-nT-r_2)S_n(r_1)S_1(r_2)x\,dr_1\,dr_2 \\ \\
\displaystyle -  \int_0^{s-nT}\int_0^{nT}(c')^+(nT-r_1,s-nT-r_2)S_n(r_1)S_1(r_2)x\,dr_1\,dr_2 \biggr) \,ds. 
\end{array}$$}
We take the second term and apply the operator $A,$ to obtain, using change of variables and Fubini's theorem, that
\small{$$\begin{array}{l}
\displaystyle A\int_0^{nT}S_n(r_1)\int_{nT}^{t}(a*1)(t-s)\int_0^{s-nT}(c')^-(nT-r_1,s-nT-r_2)S_1(r_2)x\,dr_2\,ds\,dr_1 \\ \\
\displaystyle =A\int_0^{nT}S_n(r_1)\int_{nT}^{t}(a*1)(t-s)\int_{nT}^s (c')^-(nT-r_1,u-nT)S_1(s-u)x\,du\,ds\,dr_1 \\ \\
\displaystyle =A\int_0^{nT}S_n(r_1)\int_{nT}^{t}(c')^-(nT-r_1,u-nT)\int_{u}^t (a*1)(t-s)S_1(s-u)x\,ds\,du\,dr_1 \\ \\
\displaystyle =A\int_0^{nT}S_n(r_1)\int_{nT}^{t}(c')^-(nT-r_1,u-nT)\int_{0}^{t-u} (a*1)(t-u-v)S_1(v)x\,dv\,du\,dr_1 \\ \\
\displaystyle =\int_0^{nT}S_n(r_1)\int_{nT}^{t}(c')^-(nT-r_1,u-nT)(S_1(t-u)-a(t-u))x\,du\,dr_1 \\ \\
\displaystyle =\int_0^{nT}S_n(r_1)\int_{0}^{t-nT}(c')^-(nT-r_1,t-nT-r_2)(S_1(r_2)-a(r_2))x\,dr_2\,dr_1 \\ \\
\displaystyle =((c')^-*_2 (S_n\otimes S_1))(nT,t-nT)x \\ \\
\displaystyle - \int_0^{nT}\int_{0}^{t-nT}(c')^-(nT-r_1,t-nT-r_2)a(r_2)S_n(r_1)x\,dr_2\,dr_1,
\end{array}$$}
where we have used that { $\{S_1(t)\}_{t\in (0,T]}$} is a local $(a*1,a)$-regularized resolvent family generated by $A.$ Now, observe that because { $\{S_n(t)\}_{t\in (0,nT]}$} is a local $(a*1,a)$-regularized resolvent family generated by $A,$ then $S_n(r_1)x=A(a*1*S_n)(r_1)x+a(r_1)x,$ and \small{$$\begin{array}{l}\displaystyle -\int_0^{nT}\int_{0}^{t-nT}(c')^-(nT-r_1,t-nT-r_2)a(r_2)S_n(r_1)x\,dr_2\,dr_1 \\ \\
\displaystyle=-\int_0^{nT}\int_{0}^{t-nT}(c')^-(nT-r_1,t-nT-r_2)a(r_2)(A(a*1*S_n)(r_1)+a(r_1))x\,dr_2\,dr_1.
\end{array}$$}
On the one hand, \small{$$\begin{array}{l}
\displaystyle -A\int_0^{nT}\int_{0}^{t-nT}(c')^-(nT-r_1,t-nT-r_2)a(r_2)(a*1*S_n)(r_1)x\,dr_2\,dr_1 \\ \\
\displaystyle=-A\int_0^{nT}(1*S_n)(u)x\int_{0}^{nT-u}\int_{0}^{t-nT}(c')^-(nT-u-v,t-nT-r_2)a(v)a(r_2)\,dr_2\,dv\,du \\ \\
\displaystyle=-A\int_0^{nT}((c')^- *_2 (a\otimes a))(nT-u,t-nT)(1*S_n)(u)x \,du \\ \\
\displaystyle=-A\int_0^{nT}a^-(nT-u,t-nT)(1*S_n)(u)x \,du \\ \\
\displaystyle=-A\biggl(\int_0^{2nT-t}a(2nT-t-u)(1*S_n)(u)x \,du + \int_{2nT-t}^{nT}a(t+u-2nT)(1*S_n)(u)x \,du\biggr),
\end{array}$$} and on the other hand \small{\begin{eqnarray*}
&&-\int_0^{nT}\int_{0}^{t-nT}(c')^-(nT-r_1,t-nT-r_2)a(r_1)a(r_2)x\,dr_2\,dr_1 \\ \\
&&=-((c')^- *_2 (a\otimes a))(nT,t-nT)x=-a^-(nT,t-nT)x \\ \\
&&=-a(2nT-t)x=A(a*1*S_n)(2nT-t)x-S_n(2T-t)x,
\end{eqnarray*}}
where we have applied Theorem \ref{le3.4} (ii) and that $\{S_n(t)\}_{t\in (0,nT]}$ is a local $(a*1,a)$-regularized resolvent family generated by $A.$ Then the second term is equal to \small{$$\begin{array}{l}
\displaystyle A\int_0^{nT}S_n(r_1)\int_{nT}^{t}(a*1)(t-s)\int_0^{s-nT}(c')^-(nT-r_1,s-nT-r_2)S_1(r_2)x\,dr_2\,ds\,dr_1 \\ \\
=((c')^-*_2 (S_n\otimes S_1))(nT,t-nT)x-S_n(2T-t)x \\ \\
\displaystyle -A\int_{2nT-t}^{nT}a(t+u-2nT)(1*S_n)(u)x \,du.
\end{array}$$}
Similarly, we repeat the process for the third term, and we get \small{$$\begin{array}{l}
\displaystyle -A\int_0^{nT}S_n(r_1)\int_{nT}^{t}(a*1)(t-s)\int_0^{s-nT}(c')^+(nT-r_1,s-nT-r_2)S_1(r_2)x\,dr_2\,ds\,dr_1 \\ \\
\displaystyle=-((c')^+*_2 (S_n\otimes S_1))(nT,t-nT)x-a(t)x-A\int_0^{nT}a(t-u)(1*S_n)(u)x \,du,
\end{array}$$}
where we have applied Theorem \ref{le3.4} (i). We join all the terms and we get
\small{$$\begin{array}{l}
\displaystyle A(a*1*S_{n+1})(t)x=S_{n+1}(t)x -a(t)x + A\biggl(\int_0^{nT}(a*1)(t-s)S_n(s)x\,ds \\ \\
\displaystyle -\int_{nT}^t(a*1)(t-s)S_n(2T-s)x\,ds-\int_{2nT-t}^{nT}a(t+u-2nT)(1*S_n)(u)x \,du \\ \\
\displaystyle -\int_0^{nT}a(t-u)(1*S_n)(u)x\,du \biggr).
\end{array}$$}
Note that \small{\begin{eqnarray*}
&&\int_0^{nT}(a*1)(t-s)S_n(s)x\,ds-\int_0^{nT}a(t-u)(1*S_n)(u)x\,du \\ \\
&&=\int_{t-nT}^{t}(a*1)(u)S_n(t-u)x\,du-\int_0^{nT}a(t-u)(1*S_n)(u)x\,du \\ \\
&&=(1*a)(t-nT)(1*S_n)(nT)x \\ \\
&&=\int_{2nT-t}^{nT}a(t+u-2nT)(1*S_n)(u)x \,du-\int_{0}^{t-nT}(a*1)(u)S_n(2T-t+u)x\,du \\ \\
&&=\int_{2nT-t}^{nT}a(t+u-2nT)(1*S_n)(u)x \,du-\int_{nT}^{t}(a*1)(t-s)S_n(2T-s)x\,ds,
\end{eqnarray*}}
where we have used { change of variables} and \cite[Lemma 2.2]{Mi-Po14} (note that this Lemma is true when one of the two functions is a vector valued function). Then we conclude the result.
\edproof

\section{Algebraic time translation identities for $(a,k)$-regularized resolvent families}
\label{algebraic}

\setcounter{theorem}{0}
\setcounter{equation}{0}

In this section, applying Laplace transform methods, we solve the problem of time translation in case of global $(a,k)$-regularized resolvent families. Under certain conditions on the kernels $(a,k),$ we  know that Definition \ref{ak} of $(a,k)$-regularized resolvent families is equivalent to the existence of a commutative and strongly continuous family of bounded and linear operators that satisfy $ \displaystyle \lim_{t \to 0^+} \frac{S(t)x}{k(t)} =x$ for all $x\in X$ and
the functional equation
\begin{equation}\label{eq1.1aaa}
\begin{array}{lll}
 S(s)\displaystyle \int_0^t a(t-\tau)S(\tau)x\,d\tau &-& \displaystyle  S(t) \int_0^s a(s-\tau)S(\tau)x\,d\tau   \\ \\ =k(s)\displaystyle \int_0^t a(t-\tau)S(\tau)x\,d\tau &-&k(t)\displaystyle \int_0^s a(s-\tau)S(\tau)x\,d\tau,
\end{array}
\end{equation}
for $t,s \in (0,\tau), $  see \cite[Theorem 3.1]{Liz-Po12}.

Let $S(t)$ be an $(a,k)$-regularized resolvent family in a Banach space $X$. In what follows, we will suppose that the commutative and locally integrable family $\{S(t)\}_{t > 0}$ as well as the kernels $a,k\in L^1_{loc}(\R^+)$ are Laplace transformable, with $k\in\mathcal{C}(0,\infty).$ We note that an application of the double Laplace transform to (\ref{eq1.1aaa}) gives the following identity which appears in \cite[Remark 3.2]{Liz-Po12}:
\begin{equation}\label{eq2.5}
\hat S(\lambda) \hat S(\mu)x = \frac{\hat k(\lambda)}{\hat
a(\lambda)}\displaystyle \frac{1}{\displaystyle\frac{1}{\hat
a(\lambda)}- \displaystyle \frac{1}{\hat a(\mu)}} \hat S(\mu)x -
\frac{\hat k(\mu)}{\hat a(\mu)}\displaystyle
\frac{1}{\displaystyle\frac{1}{\hat a(\lambda)}- \displaystyle
\frac{1}{\hat a(\mu)}} \hat S(\lambda)x,
\end{equation}
{ valid for} all sufficiently large $\Re \mu, \Re \lambda,$ and $x\in X.$ Using the notations in the preceding section, and the above identity, we arrive to the following notable characterization.

\begin{theorem}\label{th4.1}
A Laplace transformable and strongly continuous family of bounded and linear operators $\{S(t)\}_{t> 0}$ is an $(a,k)$-regularized resolvent family if and only if $\displaystyle \lim_{t\to 0^+}\frac{S(t)x}{k(t)}=x$ for all $x\in X$ and the following functional equation holds
\begin{equation}\label{eq4.1}
(a^+*_2 (S\otimes S))(t,s)x= k*(a*S)_t(s)x -k_t*(a*S)(s)x, \quad t,s > 0,\,x\in X.
\end{equation}
\end{theorem}

\bgproof
From (\ref{eq2.5}) we get the equivalent identity
$$\displaystyle (\hat a(\mu) -\hat a(\lambda))\hat S(\lambda) \hat S(\mu)x = \hat k(\lambda)[\hat a(\mu) \hat S(\mu)x -\hat a(\lambda) \hat S(\lambda)x] + \hat a(\lambda)[\hat k(\lambda) -\hat k(\mu)]\hat S(\lambda)x$$ valid for all $\Re \lambda ,\Re \mu $ sufficiently large. In turn, the above identity is equivalent to
$$
\begin{array}{lcl}
\displaystyle \frac{1}{\lambda -\mu}(\hat a(\mu) -\hat a(\lambda))\hat S(\lambda) \hat S(\mu)x &=& \displaystyle \frac{1}{\lambda -\mu}\hat k(\lambda)[\widehat{(a*S)}(\mu)x - \widehat{(a*S)}(\lambda)x] \\ \\ &+& \displaystyle \frac{1}{\lambda -\mu}[\hat k(\lambda) -\hat k(\mu)]\widehat{(a*S)}(\lambda)x.
\end{array}
$$
Using the identities (\ref{eq2.3}) and (\ref{eq2.4}), Proposition \ref{th3.1} (i) and uniqueness of the Laplace transform, we have the result.
\edproof

An interesting particular case is the following corollary, that we quote here for further reference.

\begin{corollary} A Laplace transformable and strongly continuous family of bounded and linear operators $\{S(t)\}_{t> 0}$ is an $(a,a)$-regularized resolvent family if and only if $\displaystyle \lim_{t\to 0^+}\frac{S(t)x}{a(t)}=x$ for all $x\in X$ and the following functional equation holds
$$\displaystyle \int_0^t \int_0^s a(t+s-r_1-r_2)S(r_1)S(r_2)x\,dr_2\, dr_2 =\displaystyle \int_0^t \int_0^s a(r_1)a(r_2)S(t+s-r_1-r_2)x\,dr_1\, dr_2 $$ for all $t,s> 0$ and $x\in X.$
\end{corollary}

\bgproof We use Corollary \ref{cor2.3} in Theorem \ref{th4.1} and the result is obtained directly.
\edproof

Our next results have the objective of extend and recover some of the results mentioned in the introduction. We will see that in order to do that, we need to impose regularity conditions on the kernels $a$ and $k$, and therefore, the results are less general than our Theorem \ref{th4.1} above.

\begin{theorem}\label{th4.6} Let $a,k \in L^1_{loc}(\mathbb{R}_+)$ be given, with $k\in\mathcal{C}(0,\infty).$ Suppose { there exist} functions $b,c\in L^1_{loc}(\R_+)$ Laplace transformable such that $c$ is absolutely { continuous on $(0,\infty)$} and $(c')^+$ is 2-Laplace transformable,  satisfying
$$
(a*b)(t)=k(t),\ (a*c)(t)=1, \qquad t>0.
$$
A Laplace transformable and strongly continuous family of bounded and linear operators $\{S(t)\}_{t> 0}$ is an $(a,k)$-regularized resolvent family if and only if $\displaystyle \lim_{t\to 0^+}\frac{S(t)x}{k(t)}=x$ for all $x\in X$ and the following functional equation holds
\begin{equation*}
( (c')^+ *_2 (S\otimes S))(t,s)x=b_t*S(s)x- b*S_t(s)x, \qquad t, s> 0,\, x\in X.
\end{equation*}
\end{theorem}

\bgproof From (\ref{eq2.5})  we obtain the equivalent identity
$$
\begin{array}{l}
\displaystyle \frac{1}{\mu -\lambda}(\lambda\widehat{c}(\lambda)-\mu\widehat{c}(\mu))\hat S(\lambda)\hat S(\mu)x = \displaystyle \frac{1}{\mu -\lambda} \widehat{b}(\lambda) \hat S(\mu)x - \frac{1}{\mu -\lambda} \widehat{b}(\mu) \hat S(\lambda)x \\ = \displaystyle \frac{1}{\mu -\lambda} \widehat{b}(\lambda)[ \hat S(\mu)x- \hat S(\lambda)x] + \frac{1}{\mu -\lambda} [\widehat{b}(\lambda) - \widehat{b}(\mu)] \hat S(\lambda)x.
\end{array}
$$
Hence, the result follows from Corollary \ref{cor3.2} (i) and formulas (\ref{eq2.3}) and (\ref{eq2.4}).
\edproof

\begin{example}{\rm
We set $a:=g_\alpha$ for $\, 0<\alpha <1$ and $k(t):= \int_0^t K(s)ds,$ { for} $t>0$. In this case  we can choose $c=g_{1-{\alpha}}$ and $b= g_{1-\alpha}*K$ satisfying the hypothesis. Therefore, we recover
the functional equation
\begin{equation}\label{eq1.6}
  \displaystyle \left(\int_t^{t+s}- \int_0^{s}\right) (g_{1-\alpha}*K)(t+s-\sigma) S(\sigma)x\, d\sigma =\displaystyle  \alpha \int_0^t \int_0^s \frac{S(r_1)S(r_2)x}{(t+s-r_1-r_2)^{1+\alpha}}\,dr_1\, dr_2,
  \end{equation}
 which appeared  in  \cite[Theorem 8]{Mei-Peng1}, and mentioned in the introduction. Here, we include as particular case the identity
  \begin{equation}\label{eq1.2}
 \left( \int_t^{t+s}  - \int_0^s\right) \frac{S(\tau)x}{(t+s-\tau)^{\alpha}}\, d\tau = \alpha \int_0^t \int_0^s \frac{S(r_1)S(r_2)x}{(t+s-r_1-r_2)^{1+\alpha}}\,dr_1\, dr_2.
  \end{equation}
 see \cite[Definition 3]{Peng-Li-12} and the more general identity developed in \cite [Theorem 5]{Li-Sun1}.}
\end{example}

\begin{example}\label{ex4.9}{\rm
Let $a:=g_{\alpha}$ and $k:=g_{\beta +1}$ where $\alpha>0, \beta >-1.$  This choosing of the pair $(a,k)$ produces the  theory of $(\alpha,\beta)$-ROF families introduced in \cite{Ch-Li10} (see also \cite[Example 3.10]{Liz-Po12} for a more general approach in terms of $(a,k)$-regularized resolvent families). As mentioned in the introduction, a time translation formula for $(\alpha, \beta)$-ROF families was developed recently { in \cite{Li-Sun1}.} Now observe that for $a$ and $k$ as before, we can choose $c=g_{1-\alpha}$ whenever $0<\alpha<1,$ and $b=g_{\beta -\alpha +1}$ whenever $\beta -\alpha >-1$ to obtain
$$
a*c= g_{\alpha}*g_{1-\alpha}=  1 \mbox{ and } a*b= g_{\alpha}*g_{\beta -\alpha +1} = g_{\beta +1}.
$$
Therefore, the hypothesis of Theorem \ref{th4.6} are satisfied  and, in consequence,  we recover the formula
\begin{eqnarray}
&& \label{eq1.5}\displaystyle \left( \int_t^{t+s}  \displaystyle  \int_0^s \right) (s+t-r)^{\beta-\alpha}S(r)x\,dr \\
&&\nonumber = \displaystyle  \alpha \frac{\Gamma(\beta -\alpha+1)}{\Gamma(1-\alpha)} \int_0^t \int_0^s \frac{S(r_1)S(r_2)x}{(t+s-r_1-r_2)^{1+\alpha}}\,dr_1\, dr_2,
  \end{eqnarray}
  whenever $\alpha \in \mathbb{R}_+ \setminus \mathbb{N}_0, \, \beta \in \mathbb{R}_+$ and $\beta -\alpha >-1$ discovered in \cite[Theorem 5]{Li-Sun1}, but now under the restrictions  $0<\alpha <1$ and $\beta-\alpha >-1$. We observe that our result correct the above formula where a more relaxed condition on $\alpha$ is assumed, namely $\alpha \in \mathbb{R}_+ \setminus \mathbb{N}_0.$ However, we note that for $ \alpha \geq 1$ the double integral on the  right hand side of (\ref{eq1.5}) diverges, as can be easily seen.}
\end{example}

Our next result widely extends the well known  semigroup functional equation to a more general class of strongly continuous families of operators. They are connected with integral equations of Volterra type, as we will see in the next section.

\begin{theorem}\label{th4.7} Let $a \in L^1_{loc}(\mathbb{R}_+)$ be given, with $a\in\mathcal{C}(0,\infty).$ Suppose there exists $c\in L^1_{loc}(\R_+)$ Laplace transformable function such that $c$ is absolutely continuous {on $(0,\infty)$},  and $(c')^+$ is 2-Laplace transformable, satisfying $(a*c)(t)=1$ for all $t> 0.$ A Laplace transformable and strongly continuous family of bounded and linear operators $\{S(t)\}_{t> 0}$ is an $(a,a)$-regularized resolvent family if and only if $\displaystyle \lim_{t\to 0^+}\frac{S(t)x}{a(t)}=x$ for all $x\in X$ and the following functional equation holds
\begin{equation*}
S(t+s)x=-((c')^+ *_2 (S\otimes S))(t,s)x, \qquad t,s> 0,\, x\in X.
\end{equation*}

\end{theorem}

\bgproof From (\ref{eq2.5})  we get the equivalent identity
$$
\frac{1}{\lambda -\mu}(\lambda\widehat{c}(\lambda)-\mu\widehat{c}(\mu))\hat S(\lambda)\hat S(\mu)x=\frac{1}{\lambda -\mu}(\hat{S}(\mu)x-\hat{S}(\lambda)x).
$$
Hence, the result follows from Corollary \ref{cor3.2} (i) and formula (\ref{eq2.2}).
\edproof

\begin{example}{\rm
If $a=g_{\alpha},$ then $\hat a(\lambda)= \frac{1}{\lambda^{\alpha}}$ and hence we can choose $c=g_{1-{\alpha}}$ with $0<\alpha<1$ satisfying the hypothesis. Note that this example recovers the functional equation
\begin{equation}\label{eq1.3}
  S(t+s)x= \frac{\alpha}{\Gamma(1-\alpha)}\int_0^t \int_0^s \frac{S(r_1)S(r_2)x}{(t+s-r_1-r_2)^{1+\alpha}}\,dr_1\, dr_2,
  \end{equation}
  for $0<\alpha<1$ as stated in \cite[Definition 2.1 (ii)]{Mei-Pe-Zh13}.}
\end{example}

Now we consider more relaxed assumptions on the functions $a$ and $k$ than the previous theorems. In contrast, the obtained functional equations are more involved and difficult to handle. The advantage is that it permits to extend the range from $0<\alpha<1$ to $1<\alpha<2$ extending recent results and producing new formulas in cases where they were not known.

\begin{theorem}\label{th4.15} Let $a,k \in L^1_{loc}(\mathbb{R}_+)$ be given, with $k\in\mathcal{C}(0,\infty).$ Suppose there exist $b,c\in L^1_{loc}(\mathbb{R}_+)$ Laplace transformable functions satisfying $(a*c)(t)=t$ and $(a*b)(t)= (1*k)(t)$ for all $t> 0.$ A Laplace transformable and strongly continuous family of bounded and linear operators $\{S(t)\}_{t> 0}$ is an $(a,k)$-regularized resolvent family if and only if $\displaystyle \lim_{t\to 0^+}\frac{S(t)x}{k(t)}=x$ for all $x\in X$ and the following functional equation holds
\begin{eqnarray*}
&&b*(1*S)_t(s)x-b_t*(1*S)(s)x=(c*S)(t)(1*S)(s)x \\ \\
&&+(1*S)(t)(c*S)(s)x-(c^+*_2 (S\otimes S))(t,s)x,
\end{eqnarray*}
for $t,s> 0,$ and $x\in X.$
\end{theorem}

\bgproof From (\ref{eq2.5})  we obtain the equivalent identity
$$
\frac{\hat{k}(\lambda)}{\hat{a}(\lambda)} \hat S(\mu)x - \frac{\hat{k}(\mu)}{\hat{a}(\mu)} \hat S(\lambda)x=(\frac{1}{\hat{a}(\lambda)}-\frac{1}{\hat{a}(\mu)})\hat S(\lambda)\hat S(\mu)x.
$$
We multiply the identity by $\frac{1}{(\lambda-\mu)\lambda \mu},$ and we obtain that \begin{eqnarray*}
\frac{1}{\lambda-\mu}\biggl(\frac{\hat{b}(\lambda)}{\mu} \hat S(\mu)x - \frac{\hat{b}(\mu)}{\lambda} \hat S(\lambda)x\biggr)&=&\biggl(\frac{\hat{c}(\lambda)}{\mu}+\frac{\hat{c}(\mu)}{\lambda}\biggr)\hat S(\lambda)\hat S(\mu)x \\
&-&\frac{1}{\lambda-\mu}\biggl( \hat{c}(\mu)-\hat{c}(\lambda) \biggr) \hat S(\lambda)\hat S(\mu)x.
\end{eqnarray*}
Note that $$\begin{array}{l}
\displaystyle\frac{1}{\lambda-\mu}\biggl(\frac{\hat{b}(\lambda)}{\mu} \hat S(\mu)x - \frac{\hat{b}(\mu)}{\lambda} \hat S(\lambda)x\biggr) \\ \\
\displaystyle =\frac{1}{\lambda-\mu}\biggl( \hat{b}(\lambda)(\widehat{(1*S)}(\mu)x -\widehat{(1*S)}(\lambda)x)- (\hat{b}(\mu)-\hat{b}(\lambda))\widehat{(1*S)}(\lambda)x \biggr).\end{array}$$
The result follows from Proposition  \ref{th3.1} (i) and formulas \eqref{eq2.3} and \eqref{eq2.4}.
\edproof

\begin{example}{\rm
Set $a=g_{\alpha}, $ with $1<\alpha <2$ and $k=g_1.$
In this case we can choose $b=c=g_{2-\alpha}$ to satisfy  $a*c= g_2$ and $a*b=1*k.$ Therefore we recover
the functional equation
\begin{eqnarray}
 \nonumber&&\displaystyle \left(\int_t^{t+s}-\int_0^{s}\right) \int_0^{\sigma} \frac{S(\tau)x}{(t+s-\sigma)^{\alpha-1}}\, d\tau\, d\sigma=\displaystyle{ \int_0^t \int_0^s \frac{S(\sigma)S(\tau)x}{(t-\sigma)^{\alpha-1}}\, d\tau\, d\sigma } \\
\label{eq1.7}&&+\displaystyle \int_0^t \int_0^s \frac{S(\sigma)S(\tau)x}{(s-\tau)^{\alpha-1}}\, d\tau\, d\sigma\displaystyle{- \int_0^t \int_0^s \frac{S(\sigma)S(\tau)x}{(t+s-\sigma-\tau)^{\alpha-1}}\, d\tau\, d\sigma}.
 \end{eqnarray}
developed  in \cite[Definition 3.1]{Mei-Peng-Jia} and cited in the introduction.} \end{example}

\begin{example}{\rm
Let $a=g_{\alpha}$ and $ k=g_{\beta +1}$ where $\alpha>0, \beta >-1.$ In this case we obtain a new functional equation for $(\alpha,\beta)$-ROF families (introduced in \cite{Ch-Li10}) in contrast with those developed in \cite{Li-Sun1}. See also Example \ref{ex4.9} for a correction on the assumptions on $\alpha$ and $\beta$. Indeed, we can choose $c=g_{2-\alpha}$ and $b=g_{\beta -\alpha +2}$ under the assumptions $0<\alpha<2$ and $\beta-\alpha>-2.$
}
\end{example}

\section{Examples,  applications and final comments.}

\setcounter{theorem}{0} \setcounter{equation}{0}

\subsection{Multiplication local regularized families in $L^p(\R)$}

We consider the Lebesgue space $L^p(\R),$ $1\leq p\leq \infty,$ and $a=g_{\alpha}$ with $\alpha\in (0,2).$ Define the multiplication operator $$Af(x):=(1+x+ix^2)^{\alpha},\qquad x\in\R,\, f\in L^p(\R),$$ with maximal domain in $L^p(\R).$  Assume $s\in(1,2),$ $\delta=\frac{1}{s}$ and $K_{\delta}(t)=\mathcal{L}^{-1}(e^{-\lambda^{\delta}})(t),$ \mbox{$t> 0,$} where $\mathcal{L}^{-1}$ is the inverse Laplace transform. { Then $A$ generates a global $(a,K_{\delta})$-regularized resolvent family in $L^p(\R).$ Furthermore, when $s=2$ there exists $\tau>0$ such that $A$ generates a local $(a,K_{\frac{1}{2}})$-regularized resolvent family on $[0,\tau),$ see \cite[Example 2.31]{Kostic}. Then we can apply Theorem \ref{local1}, and conclude that $A$ generates a local $(a,(K_{\frac{1}{2}}*a)^{*n}*K_{\frac{1}{2}})$-regularized resolvent family on $(0,(n+1)\tau),$ for all $n\in\mathbb{N}.$}

\subsection{Local regularized families in { sequence} spaces}

Let $l^2(\NN)=\{ x=(x_m)_{m=1}^{\infty}\subset\CC\,:\, \displaystyle\sum_{n=1}^{\infty}|x_m|^2<\infty \}$ be the Hilbert space of all square-summable sequences with the norm $\lVert x\rVert=(\displaystyle\sum_{n=1}^{\infty}|x_m|^2)^{1\over 2}.$ We take $\tau>0,$ and $$a_m=\frac{m}{\tau}+i\biggl(\biggl( \frac{e^m}{m}\biggr)^2-\biggl(\frac{m}{\tau}\biggr)^2 \biggr)^{\frac{1}{2}}, \qquad m\in \NN,$$ where $i$ is the imaginary identity.{ We note that for all $n\in \NN,$ the sequence $(a_m)_{m=1}^{\infty}$ generates a local $n$-times integrated semigroup on $l^2(\NN)$ for $t\in [0,n\tau),$ see \cite[p.75-76]{Oka-Ta90}.}

The Mittag-Leffler functions are defined by $$E_{\alpha,\beta}(z):=\displaystyle\sum_{n=0}^{\infty}\frac{z^n}{\Gamma(\alpha n+\beta)}, \qquad \alpha,\beta>0,z\in\mathbb{C}.$$ For short, $E_{\alpha}:=E_{\alpha,1}.$ We take $0<\alpha<2.$ Observe that the function $E_\alpha$ is a $(g_{\alpha},1)$-regularized resolvent family.

For any $\beta\in\R^+,$ $(T_{\alpha,\beta}(t))_{t\in (0,\beta \tau)},$ defined by $$T_{\alpha,\beta}(t)x=\biggl( \frac{1}{\Gamma(\beta)}\int_0^t(t-s)^{\beta-1}E_{\alpha}((a_m s)^{\alpha})x_m\,ds \biggr)_{m=1}^{\infty},\text{ for }x \in l^2(\NN),$$ is a local $(g_{\alpha},g_{\beta+1})$-regularized resolvent family on $l^2(\NN):$

Note that $a_m\in\CC^{+},$ the set of imaginary numbers with positive real part. Then for all $s\geq 0$ and  $0<\alpha<2,$ $|arg\,(a_m s)^{\alpha}|\leq \frac{\alpha \pi}{2}.$ So, the asymptotic expansion \cite[(1.27)]{Bajlekova} and the continuity of the Mittag-Leffler function imply that there are constants $c,C$ such that $c e^{a_m s} \leq E_{\alpha}((a_m s)^{\alpha})\leq C e^{a_m s}.$ Observe that $$\lVert T_{\alpha,\beta}(t)\rVert=\displaystyle\sup_{m\in\NN}|\frac{1}{\Gamma(\beta)}\int_0^t(t-s)^{\beta-1}E_{\alpha}((a_m s)^{\alpha})\,ds|<\infty$$ if and only if $$\displaystyle\sup_{m\in\NN}|\frac{1}{\Gamma(\beta)}\int_0^t(t-s)^{\beta-1}e^{a_m s}\,ds|<\infty,$$ which happens if only if $0\leq t< \beta \tau,$ see \cite[Example 1]{Miana}. It is clear that $\{T_{\alpha,\beta}(t) \}_{t\in (0,\beta \tau)}$ is strongly continuous and verifies any functional equation associated to $(g_{\alpha},g_{\beta+1})$-regularized families. The case $\alpha=1$ is made in \cite{Miana}.

\subsection{A new class of regularized families without jumps of  regularity}

Let $0<\tau\leq \infty$ and $b\in\ L^1_{loc}(\R^+).$ We define $a=b*b.$ Suppose that $\{R_1(t)\}_{t\in (0,\tau)}$ is a local $(a,a)$-regularized resolvent family generated by $A,$ such that $A$ verifies condition (H5) of \cite{Kostic} (for example $A$ densely defined). Then by \cite[Theorem 2.34]{Kostic}, we have that \begin{displaymath}
\mathscr{A}\equiv
\left( \begin{array}{cc}
0 & I \\
A & 0
\end{array} \right)
\end{displaymath}
is the generator of a local $(b,b^{*3})$-regularized resolvent family $\{S_1(t)\}_{t\in (0,\tau)}$ given explicitly by \begin{displaymath}
S_1(t)=
\left( \begin{array}{cc}
(b*R_1)(t) & (a*R_1)(t) \\
R_1(t)-a(t)I & (b*R_1)(t)
\end{array} \right), \qquad 0<t<\tau.
\end{displaymath}

Now, we suppose that there is a $c\in L^1_{loc}(\R^+)$ Laplace transformable such that $c$ is absolutely continuous, differentiable a.e., and $(c')^+$ is 2-Laplace transformable, and satisfying $(a*c)(t) = 1$ for all $t>0 .$ Then we conclude $A$ generates a global $(a,a)$-regularized resolvent family $\{R(t)\}_{t\in (0,\infty)}$ which extends $\{R_{1}(t)\}_{t\in (0,\tau)},$ see Theorem \ref{th-local1}. Then, we can extend $\{S_1(t)\}_{t\in (0,\tau)}$ without loss of regularity, i.e, $\mathscr{A}$ generates a global $(b,b^{*3})$-regularized resolvent family $\{S(t)\}_{t\in (0,\infty)}$ given by \begin{displaymath}
S(t)=
\left( \begin{array}{cc}
(b*R(t) & (a*R)(t) \\
R(t)-a(t)I & (b*R)(t)
\end{array} \right), \qquad 0< t.
\end{displaymath}
In the particular case $b=g_{\frac{\alpha}{2}}$ with $0<\alpha<1,$ $a=g_{\alpha}$ and $\{R_1(t)\}_{t\in (0,\tau)}$ be a local $(g_{\alpha},g_{\alpha})$-regularized resolvent family generated by $A,$ we can extend $\{S_1(t)\}_{t\in (0,\tau)}$ without loss of regularity, i.e, $\mathscr{A}$ generates a global $(g_{\frac{\alpha}{2}},g_{\frac{3\alpha}{2}})$-regularized resolvent family $\{S(t)\}_{t>0}$ such that $S(t)=S_1(t)$ for $0< t<\tau.$

\subsection{Applications to obtain new functional equations}

We give several examples of the abstract results in section 6. They show that we can recover, extend and produce new functional equations that in some cases are interesting for their own nature.

\begin{example}
{\bf $C_0$-semigroups.} \rm Choose $a(t)= 1$ and $k(t)= 1$ for $t>0$. Then, for $x\in X,$ we obtain
$$
\int_0^t \int_0^s S(\tau_1)S(\tau_2)x\, d\tau_2\, d\tau_2 = \int_t^{t+s} \int_0^r S(\tau)x\,d\tau\, dr - \int_0^s \int_0^r S(\tau)x\, d\tau\, dr, \qquad t,s>0.
$$
Taking the derivative one time with respect to the variable $t$ we obtain $$ \displaystyle S(t)\int_0^s S(\tau_2)x\,d\tau_2 = \int_0^{t+s} S(\tau)x\, d\tau -\int_0^t S(\tau)x\,d\tau, \quad t,s>0,$$
 which was introduced in \cite[Example 3.4]{Liz-Po12}. Hence taking the derivative with respect to the variable $s$ we get the Cauchy formula $ S(s)S(t)= S(t+s).$ Observe that in this way  we deduce easily Cauchy's functional equation from the formula in the Theorem \ref{th4.1}.
\end{example}

\begin{example}
{\bf Cosine families.} {\rm Choose $a(t)=t $ and $k(t)= 1$ for $t>0$. Then, for $x\in X$ and  $t,s >0$  we have
$$
\begin{array}{l}
\displaystyle \int_0^t \int_0^s (t+s-r_1-r_2)S(r_1)S(r_2)x\,dr_2\, dr_1
\displaystyle = \left(\int_t^{t+s}  - \int_0^s \right)\int_0^r (r-\tau)S(\tau)x\,d\tau dr.
\end{array}
$$
 To see directly why the above formula is equivalent to the D'Alembert functional equation { $ S(t+s) +S(|t-s|)= 2S(t)S(s)$} we proceed as in the above example, first taking derivative with respect to the variable $t,$ to obtain
$$
\begin{array}{l}
\displaystyle S(t)\int_0^s(s-r_2)S(r_2)x\,dr_2 + \int_0^t  S(r_1)\int_0^s S(r_2)x\, dr_2dr_1 \\ \\
\displaystyle = \int_0^{t+s}(t+s-\tau)S(\tau)x\,d\tau -\int_0^t(t-\tau)S(\tau)x\,d\tau,
\end{array}
$$
and then taking derivative with respect to the variable $s,$ to have
$$
S(t)\int_0^s S(r_2)x\,dr_2 + \int_0^t S(r_1)S(s)x\,dr_1 = \int_0^{t+s}S(\tau)x\,d\tau, \qquad t,s>0.
$$
By  \cite[Theorem 2]{Mei-Pe-Zh(SemFor)} { the above} functional equation is equivalent to the cosine functional equation.}
\end{example}

\begin{example}{\bf Convoluted semigroups.} {\rm
Choosing  $a(t)=1$ for $t>0$ and $k \in C^2(\mathbb{R}_+)$ and proceeding as in the above examples, we obtain a new functional equation for convoluted semigroups:
\small{$$\begin{array}{lcl}
S(t)S(s)x &=& k(0)S(t+s)x +\displaystyle k'(0)\int_0^{t+s} S(\tau)x\, d\tau -k'(s)\int_0^t S(\tau)x\,d\tau \\ \\ &-& \displaystyle k'(t)\int_0^s S(\tau)x\, d\tau + \int_t^{t+s}k''(t+s-r)\int_0^r S(\tau)x\,d\tau\, dr \\ \\ &-& \displaystyle \int_0^s k''(t+s-r)\int_0^r S(\tau)x\,d\tau\, dr,
\end{array}$$}
for $t,s>0$ and $x\in X.$ Setting $k(t):=(1-\epsilon) + \epsilon t $ for $0\leq \epsilon \leq 1,$  and $t\ge 0$, this formula shows an interesting fact: How  continuously moves the functional equation from the case of semigroups to the case of $1$-times integrated semigroup as $\epsilon$ varies from $0$ to $1$:
$$
S(t)S(s)x = (1-\epsilon)S(t+s)x + \epsilon\left( \left(\int_t^{t+s}  -\int_0^s \right)S(\tau)x\,d\tau\right), \qquad t,s>0, \quad x\in X.
$$
}

\end{example}

\begin{example}{\bf Resolvent families.} {\rm  Now we take $k(t)= 1$ for $t>0$; $a\in C^2(\mathbb{R_+})$ and
 we proceed as above. We get this  (new) functional equation:
\small{$$\begin{array}{lcl}
a(0)S(t)S(s)x &+& \displaystyle S(t)\int_0^s a'(s-\tau)S(\tau)x\,d\tau + S(s)\int_0^t a'(t-\tau)S(\tau)x\,d\tau \\ \\
&+& 2 \displaystyle \int_0^t \int_0^s a''(t+s-r_1-r_2)S(r_1)S(r_2)x\,dr_2\, dr_1 \\ \\
&=& a(0)S(t+s)x + \displaystyle \int_0^s a'(t+s-\tau)S(\tau)x\,d\tau,
\end{array}$$}
for $x\in X$. For $0\leq \epsilon \leq 1,$ we set $a(t):= (1-\epsilon) + \epsilon t $, and we see how the formula continuously moves from the semigroup to the cosine family cases when $\epsilon$ goes from $0$ to $1$:
$$
(1-\epsilon)[S(t)S(s)x -S(t+s)x] = \epsilon[\int_0^{t+s} S(\tau)x\,d\tau -\int_0^s S(\tau)x\,d\tau - \int_0^t S(\tau)x\,d\tau],
$$ for $x\in X.$ Note the intriguing case  $\epsilon =1/2$ where the difference between both functional equations is the same,  which indicates that  in some suitable norm the study of the topology of the set of all functional equations satisfying (\ref{eq1.1}) should be relevant.}

\end{example}

{\bf Acknowledgement.} We  thank an anonymous referee for several comments, ideas and references which have contributed to improve the final version of the paper.

\bibliographystyle{amsplain}

\end{document}